\documentclass{amsart}
\usepackage{amsmath, amssymb}
\usepackage{amsfonts}
\usepackage[arrow,matrix,curve,cmtip,ps]{xy}

\usepackage{amsthm}

\allowdisplaybreaks

\newtheorem{theorem}{Theorem}[section]
\newtheorem{lemma}[theorem]{Lemma}
\newtheorem{proposition}[theorem]{Proposition}
\newtheorem{corollary}[theorem]{Corollary}

\newtheorem*{theorem*}{Theorem}
\theoremstyle{remark}
\newtheorem{remark}[theorem]{Remark}
\newtheorem{definition}[theorem]{Definition}
\newtheorem{example}[theorem]{Example}

\numberwithin{equation}{section}

\newcommand{\W}{\text{Wojciech}}
\newcommand{\wvec}{\omega}
\newcommand{\wmap}{\omega}
\newcommand{\Z}{\mathbb{Z}}
\newcommand{\N}{\mathbb{N}}

\newcommand{\K}{\mathcal{K}}
\newcommand{\B}{\mathcal{B}}
\newcommand{\Q}{\mathcal{Q}}
\newcommand{\Hi}{\mathcal{H}}

\newcommand{\im}{\operatorname{im }}
\newcommand{\coker}{\operatorname{coker }}
\newcommand{\Ad}{\operatorname{Ad}}
\newcommand{\ext}{\operatorname{Ext}}

\newcommand{\ind}{\operatorname{ind}}
\newcommand{\rank}{\operatorname{rank}}

\begin{document}
\title{Computing Ext for Graph Algebras}
\author{Mark Tomforde 
}

\address{Department of Mathematics\\ Dartmouth College\\
Hanover\\ NH 03755-3551\\ USA}

\curraddr{Department of Mathematics\\ University of Iowa\\
Iowa City\\ IA 52242\\ USA}

\email{tomforde@math.uiowa.edu}

\thanks{This research was carried out while the author was a
student at Dartmouth College and it forms part of his doctoral
dissertation.  The author would like to take this
opportunity to thank Dana P.~Williams for his supervision and
guidance throughout this project.  The author would also like
to thank the referee for providing many helpful suggestions.}

\date{\today}
\subjclass{19K33 \and 46L55}

\begin{abstract}
\noindent
For a row-finite graph $G$ with no sinks and in which
every loop has an exit, we construct an isomorphism
between $\ext (C^*(G))$ and $\coker (A-I)$, where $A$ is the
vertex matrix of $G$.  If $c$ is the class in $\ext(C^*(G))$
associated to a graph obtained by attaching a sink to $G$, then
this isomorphism maps $c$ to the class of a vector that
describes how the sink was added.  We conclude with an
application in which we use this isomorphism to produce an
example of a row-finite transitive graph with no sinks whose
associated $C^*$-algebra is not semiprojective. 
\end{abstract}

\maketitle

\section{Introduction}

The Cuntz-Krieger algebras $\mathcal{O}_A$ are
$C^*$-algebras that are generated by a collection of
partial isometries satisfying relations described by a
finite matrix $A$ with entries in $\{ 0 , 1\}$ and no zero
rows.  In \cite{CK} Cuntz and Krieger computed $\ext$ for
these $C^*$-algebras, showing that $\ext \mathcal{O}_A$
is isomorphic to $\coker (A-I)$, where $A : \Z^n \rightarrow
\Z^n$.

In 1982 Watatani noted that one can view $\mathcal{O}_A$ as
the $C^*$-algebra of a finite directed graph $G$ with no
sinks and whose vertex adjacency matrix is $A$ \cite{Wat}. 
However, it was not until the late 1990's that  analogues of
these $C^*$-algebras were considered for possibly infinite
graphs that are allowed to contain sinks \cite{KPR, KPRR}. 
Since that time there has been a flurry of activity in
studying these graph algebras. 

Graph algebras have proven to be important for
many reasons.   To begin with, they include a fairly wide
class of $C^*$-algebras.  In addition to generalizing the
Cuntz-Krieger algebras, graph algebras include many
other interesting classes of $C^*$-algebras such as AF-algebras
and Kirchberg-Phillips algebras with free
$K_1$-group.  However, despite the fact that graph algebras
include a wide class of $C^*$-algebras, their basic structure
is fairly well understood and their invariants are readily
computable.  In fact, results about Cuntz-Krieger
algebras can often be extended to graph algebras with only minor
modifications.  One reason graph algebras have attracted
the interest of many people is that the graph provides a
convenient tool for visualization.  Not only does the graph
determine the defining relations for the generators of the
$C^*$-algebra, but also many important properties of the
$C^*$-algebra may be translated into graph properties that
can easily be read off from the graph.    

In this paper we extend Cuntz and Krieger's computation of
$\ext \mathcal{O}_A$ to graph algebras.  Specifically, we
prove the following.
\begin{theorem*}  Let $G$ be a row-finite graph with
no sinks and in which every loop has an exit, and let
$C^*(G)$ be the $C^*$-algebra associated to $G$.  Then there
exists an isomorphism $$\wmap: \ext (C^*(G))
\rightarrow
\coker (A_G-I)$$ where $A_G$ is the vertex matrix of $G$ and
$A_G : \prod_{G^0} \Z \rightarrow \prod_{G^0} \Z$. 
\end{theorem*}

In addition to showing that $\ext (C^*(G)) \cong
\coker (A_G-I)$, the isomorphism $\wmap$ is important because
its value on certain extensions can be easily calculated.  If
$E$ is an essential 1-sink extension of $G$ as described in
\cite{RTW}, then $C^*(E)$ will be an extension of $C^*(G)$ by
$\K$ and thus determines an element in $\ext(C^*(G))$. 
Roughly speaking, a 1-sink extension of
$G$ may be thought of as a graph formed by attaching a sink
$v_0$ to $G$, and this 1-sink extension is said to be
essential if every vertex of $G$ can reach this sink.  For
example, if
$G$ is the graph
\vspace{.2in}
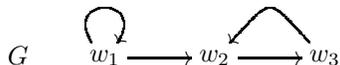
$$
\begin{matrix}
G & &
\xymatrix{
w_1 \ar@(ul,ur) \ar[r] & w_2 \ar[r] & w_3 \ar@(ul,ur)[l] \\
} \end{matrix}$$ then two examples of essential 1-sink
extensions are the following graphs $E_1$ and $E_2$.
\vspace{.2in}
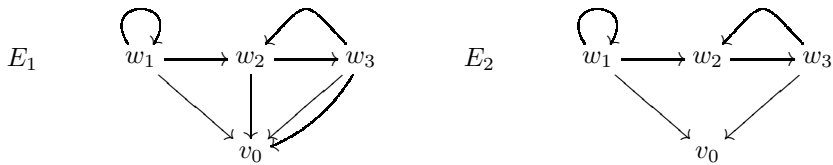
$$ \begin{matrix} E_1 & & &
\xymatrix{
w_1 \ar@(ul,ur) \ar[r] \ar[dr] & w_2 \ar[r] \ar[d] & w_3
\ar@(ul,ur)[l] \ar[dl] \ar@/^/[dl] \\
 & v_0 &  \\
}
& & & E_2 & & &
\xymatrix{
w_1 \ar@(ul,ur) \ar[r] \ar[dr]  & w_2 \ar[r]  &
w_3
\ar@(ul,ur)[l] \ar[dl]  \\
 & v_0 &  \\
}
\end{matrix}
$$
For each 1-sink extension there is a vector, called
the $\W$ vector, that describes how the sink is added
to $G$ \cite{RTW}.  In the above two examples the $\W$ vector
is the vector whose $v^{\text{th}}$ entry is equal to the
number of edges from
$v$ to the sink.  This vector is $\left(
\begin{smallmatrix} 1 \\ 1 \\ 2 \\
\end{smallmatrix} \right)$ for $E_1$ and $\left(
\begin{smallmatrix} 1 \\ 0 \\ 1 \\
\end{smallmatrix} \right)$ for $E_2$.  It turns out that if
$E$ is a 1-sink extension of $G$, then the value that $\wmap$
assigns to the element of
$\ext(C^*(G))$ associated to $E$ is equal to the class of the
$\W$ vector of $E$ in $\coker (A_G-I)$.  Furthermore, since
$\wmap$ is additive we have a nice way of describing
addition of elements in $\ext(C^*(G))$ associated to
essential 1-sink extensions.  For example, if $E_1$ and $E_2$
are as above, then the sum of their associated elements in
$\ext(C^*(G))$ is the element in $\ext(C^*(G))$ associated to
the 1-sink extension
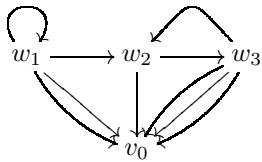
\vspace{.2in}
$$ \xymatrix{
w_1 \ar@(ul,ur) \ar[r] \ar[dr] \ar@/_/[dr] & w_2
\ar[d] \ar[r]  & w_3
\ar@(ul,ur)[l] \ar[dl] \ar@/_/[dl] \ar@/^/[dl] \\
 & v_0 &  \\
}
$$
whose $\W$ vector is $\left(
\begin{smallmatrix} 2 \\ 1 \\ 3 \\
\end{smallmatrix} \right) = \left(
\begin{smallmatrix} 1 \\ 1 \\ 2 \\
\end{smallmatrix} \right) +  \left(
\begin{smallmatrix} 1 \\ 0 \\ 1 \\
\end{smallmatrix} \right)$.  Thus we have a way of
visualizing certain elements of $\ext (C^*(G))$ as well as a
way to visualize their sums.  We show in \S \ref{Wmapsection}
that if $G$ is a finite graph, then every element of
$\ext(C^*(G))$ is an element associated to an essential 1-sink
extension of $G$.  We also show that this is not
necessarily the case for infinite graphs.

In addition to providing an easily visualized description of
$\ext (C^*(G))$, we also show that the isomorphism
$\wmap$ can be used to ascertain information about the
semiprojectivity of a graph algebra.  Blackadar has shown
that the Cuntz-Krieger algebras are semiprojective
\cite{Bla2}, and Szyma\'nski has proven that $C^*$-algebras
of transitive graphs with finitely many vertices are
semiprojective \cite{Szy}.   Although not all graph
algebras are semiprojective (for instance, it follows from
\cite[Theorem 3.1]{Bla2} that $\K$ is not semiprojective), it is
natural to wonder if the $C^*$-algebras of transitive graphs
will always be semiprojective.  In
\S \ref{semiprojectivity} we answer this question in the
negative.  We use the isomorphism $\wmap$ to produce an example
of a row-finite transitive graph whose
$C^*$-algebra is not semiprojective.

This paper is organized as follows.  We begin in
\S\ref{extprelim} with a description of
$\ext$ due to Cuntz and Krieger.  Then, after some graph algebra
preliminaries in \S\ref{graphalgprelimsec}, we continue in
\S\ref{CKmapsec} by defining a map $d : \ext (C^*(G))
\rightarrow \coker (B_G-I)$, where $B_G$ is the edge matrix of
$G$.  In \S\ref{Wmapsection} we define the map $\wmap :
\ext(C^*(G)) \rightarrow \coker (A_G-I)$, where $A_G$ is the
vertex matrix of $G$.  We also prove that $\wmap$ is an
isomorphism and compute the value it assigns to elements of
$\ext(C^*(G))$ associated to essential 1-sink extensions.  We
conclude in \S\ref{semiprojectivity} by providing an example of
a row-finite transitive graph whose $C^*$-algebra
is not semiprojective.
 
This research was carried out while the author was a student at
Dartmouth College and it forms part of his doctoral
dissertation.  The author would like to take this opportunity
to thank Dana P. Williams for his supervision and guidance
throughout this project.  The author would also like to thank
the referee for providing many helpful suggestions.

\section{$\ext$ Preliminaries}
\label{extprelim}

Throughout we shall let $\Hi$ denote a separable
infinite-dimensional Hilbert space, $\K$ denote the compact
operators on $\Hi$, $\B$ denote the bounded operators on $\Hi$
, and $\Q := \B / \K$ denote the associated Calkin algebra.  We
shall also let $i : \K \rightarrow \B$ denote the inclusion map
and $\pi : \B \rightarrow \Q$ denote the projection map.

In this section we review a few definitions and establish
notation.  We asume that the reader is familiar with $\ext$. 
For those readers who would like more background on $\ext$ we
suggest \cite{Bla} and \cite{JT}, or for a less comprehensive
but more introductory treatment we suggest \cite{WO}.  We also
mention that an expanded version of the topics addressed here,
including an account of $\ext$, is contained in \cite{thesis}.

If $A$ is a $C^*$-algebra, then an \emph{extension} of $A$ (by
the compact operators) is a homomorphism $\tau : A \rightarrow
\Q$.  An extension is said to be \emph{essential} if it is a
monomorphism.
\begin{definition}  An extension $\tau : A \rightarrow \Q$ is
said to be \emph{degenerate} if there exists a homomorphism
$\eta : A \rightarrow \B$ such that $\pi \circ \eta = \tau$.
In other words, $\tau$ can be lifted to a (possibly degenerate)
representation $\eta$.
\end{definition}
We warn the reader that the terminology used above is not
standard.  Many authors refer to such extensions as
trivial rather than degenerate.  However, we have chosen to
follow the convention established in \cite{JT}.

It is a fact that if there exists an essential degenerate
extension of $A$ by $\K$ then $\ext(A)$ will be comprised of
weakly stable equivalence classes of essential extensions
\cite[Proposition 15.6.5]{Bla}.  However, we will find it more
convenient to use a description of $\ext$ given by Cuntz and
Krieger in \cite{CK} when they computed $\ext \mathcal{O}_A$.

\begin{definition}
We say that two Busby invariants $\tau_1$ and $\tau_2$ are
\emph{CK-equivalent} if there exists a partial isometry
$v \in \Q$ such that
\begin{equation} \tau_1 = \mathrm{Ad}(v) \circ \tau_2
\hspace{.15in} \text{ and } \hspace{.15in} \tau_2 =
\mathrm{Ad}(v^*)
\circ
\tau_1.
\label{CK-equiv} \end{equation}
\label{CK-eq}
\end{definition}

The following fact was used in \cite{CK}.

\begin{lemma}  Suppose that $\tau_1$ and $\tau_2$ are the
Busby invariants of two essential extensions of $A$ by $\K$. 
Then $\tau_1$ equals $\tau_2$ in $\ext (A)$ if and
only if $\tau_1$ and $\tau_2$ are CK-equivalent.
\label{CK-equivalent}
\end{lemma}

\noindent In light of this lemma we may think of the class of
$\tau$ in $\ext (A)$ as the class generated by the
relation in (\ref{CK-equiv}).  Furthermore, we see that
any two essential degenerate extensions will be equivalent.

For extensions $\tau_1$ and $\tau_2$ we say that $\tau_1 \perp
\tau_2$ if there are orthogonal projections $p_1$ and $p_2$
such that $\tau_i (A) \subseteq p_i
\Q p_i$.  In this case we may define a map $\tau_1 \boxplus
\tau_2$ by $a \mapsto \tau_1(a) + \tau_2(a)$.  The
orthogonality of the projections is enough to ensure that
this map will be multiplicative and therefore $\tau_1 \boxplus
\tau_2$ will be a homomorphism.  The notation $\boxplus$
is used because a quite different
meaning has already been assigned to $\tau_1 + \tau_2$ in $\ext
(A)$.

Provided that there exists an essential degenerate extension of
$A$ by $\K$, we may view $\ext (A)$ as the
equivalence classes of essential extensions generated
by the relation in (\ref{CK-equiv}).  For any two elements
$\tau_1 , \tau_2  \in \ext(A)$, we define their sum
to be $ \tau_1 + \tau_2  = \tau_1' \boxplus \tau_2'$ where
$\tau_1'$ and $\tau_2'$ are essential extensions such that
$\tau_1' \perp \tau_2'$ and
$\tau_i'$ is weakly stably equivalent to $\tau_i$.  Note that
the common class of all degenerate essential extensions acts as
the neutral element in $\ext (A)$.

\section{Preliminaries on Graph $C^*$-Algebras}
\label{graphalgprelimsec}

A (directed) graph $G = (G^0,G^1,r,s)$ consists of
a countable set $G^0$ of vertices, a countable set $G^1$ of
edges, and maps $r,s:G^1 \rightarrow G^0$ that identify the
range and source of each edge.   A vertex
$v \in G^0$ is called a \emph{sink} if
$s^{-1}(v)=\emptyset$ and a \emph{source} if $r^{-1}(v) =
\emptyset$.  All of our graphs will be
assumed to be \emph{row-finite} in that each vertex emits only
finitely many edges 

If $G$ is a row-finite directed graph, a \emph{Cuntz-Krieger
$G$-family} in a $C^*$-algebra is a set of mutually orthogonal
projections $\{ p_v : v \in G^0 \}$ together with a set of
partial isometries $\{ s_e : e \in G^1 \}$ that satisfy the
\emph{Cuntz-Krieger relations} 
$$s_e^* s_e = p_{r(e)} \ \text{for} \  e \in E^1 \ \
\text{and} \
\ p_v =
\sum_{ \{ e : s(e) =v \} } s_e s_e^* \ \text{whenever } v \in G^0
\text{ is not a sink.}$$
\noindent Then $C^*(G)$ is defined to be the $C^*$-algebra
generated by a universal Cuntz-Krieger $G$-family
\cite[Theorem 1.2]{KPR}.

A \emph{path} in a graph $G$ is a finite sequence of edges
$\alpha := \alpha_1 \alpha_2 \ldots \alpha_n$ for which
$r(\alpha_i) = s(\alpha_{i+1})$ for $1 \leq i \leq n-1$, and we
say that such a path has length $|\alpha|=n$.  For $v,w \in
G^0$ we write $v \geq w$ to mean that there exists a path
with source $v$ and range $w$.  For $K,L
\subseteq G^0$ we write $K \geq L$ to mean that for each $v \in
K$ there exists $w \in L$ such that $v \geq w$.

A \emph{loop} is a path whose range and source are equal.  An
\emph{exit} for a loop $x := x_1 \ldots x_n$ is an edge $e$ for
which $s(e)=s(x_i)$ for some $i$ and $e \neq x_i$.  A graph
is said to satisfy Condition~(L) if every loop in $G$ has an
exit.

If $G$ is a graph then we may associate two matrices to $G$. 
The \emph{vertex matrix} of $G$ is the $G^0 \times G^0$ matrix
$A_G$ whose entries are given by $A_G(v,w) := \# \{ e \in G^1 :
s(e) = v \text{ and } r(e) = w \}$.  The \emph{edge matrix} of
$G$ is the $G^1 \times G^1$ matrix
$B_G$ whose entries are given by $$B_G(e,f) := \begin{cases}
1 & \text{ if $r(e) = s(f)$.} \\ 0 & \text{
otherwise. } \end{cases}$$  Notice that if $G$ is a row-finite
graph, then the rows of both $A_G$ and $B_G$ will eventually be
zero.  Hence left multiplication gives maps $A_G : \prod_{G^0}
\Z \rightarrow \prod_{G^0} \Z$ and $B_G : \prod_{G^1} \Z
\rightarrow \prod_{G^1} \Z$.  Also the maps $A_G - I:
\prod_{G^0} \Z
\rightarrow \prod_{G^0} \Z$ and $B_G -I : \prod_{G^1} \Z
\rightarrow \prod_{G^1} \Z$ will
prove important in later portions of this paper.

\section{The Ext Group for $C^*(G)$}
\label{CKmapsec}

The proofs of the following three lemmas are straightforward.

\begin{lemma}  Suppose that $p_1,p_2,\ldots $ is a countable
sequence of pairwise orthogonal projections in $\Q$. 
Then there are pairwise orthogonal projections
$P_1,P_2,\ldots $ in $\B$ such that $\pi(P_i)=p_i$ for $ i =
1, 2, \ldots$.
\label{manyprojectionlift}
\end{lemma}

\begin{lemma}  If $w$ is a partial isometry in $\Q$, then there
exists a partial isometry $V$ in $\B$ such that $\pi(V)=w$. 
\label{partialisometrylift}
\end{lemma}

\begin{lemma}  If $w$ is a unitary in $\Q$, then $w$ can be
lifted to either an isometry or coisometry $U \in \B$.
\end{lemma}

For the rest of this section let $G$ be a row-finite graph with
no sinks that satisfies Condition~(L).  Since $C^*(G)$ is
separable, there will exist an essential degenerate
extension of $C^*(G)$ \cite[\S 15.5]{Bla}.  (In fact, we shall
prove that there are many essential degenerate extensions in
Lemma \ref{texists}.)  Therefore we may use Cuntz and Krieger's
description of $\ext$ discussed in \S\ref{extprelim}.  

Let $E \in \Q$ be a projection.  By Lemma
\ref{manyprojectionlift} we know that there exists a projection
$E' \in \B$ such that
$\pi(E')=E$.  If $X$ is an element of $\Q$ such that $EXE$ is
invertible in $E \Q E$, then we denote by $\text{ind}_E(X)$ the
Fredholm index of $E'X'E'$ in $\text{im }E'$, where $X' \in \B$
is such that $\pi(X')=X$.  Since the Fredholm index is invariant
under compact perturbations, this definition does not depend on
the choice of $E'$ or $X'$.

The following two lemmas are taken from \cite{CK}.

\begin{lemma} Let $E,F \in \Q$ be orthogonal projections, and
let $X$ be an element of $\Q$ such that $EXE$ and $FXF$ are
invertible in $E \Q E$ and $F \Q F$ and such that $X$ commutes
with $E$ and $F$.  Then $\mathrm{ind}_{E+F}(X) =
\mathrm{ind}_E(X)+\mathrm{ind}_F(X)$.
\label{CKlemma1}
\end{lemma}
\begin{lemma} Let $X$ and $Y$ be invertible
operators in $E \Q E$.  Then $\mathrm{ind}_E(XY) =
\mathrm{ind}_E(X)+\mathrm{ind}_E(Y)$.
\label{CKlemma2}
\end{lemma}

In addition, we shall make use of the following lemmas to
define a map from $\ext (C^*(G))$ into $\coker (B_G-I)$. The
first lemma is an immediate consequnce of the Cuntz-Krieger
Uniqueness Theorem for graph algebras \cite[Theorem~3.1]{BPRS}.

\begin{lemma}  Let $G$ be a graph that satisfies Condition
(L), and let $\{s_e, p_v \}$ be the canonical Cuntz-Krieger
$G$-family in $C^*(G)$.  If $I$ is an ideal of $C^*(G)$ with
the property that $p_v \notin I$ for all $v \in G^0$, then $I =
\{ 0 \}$.
\label{projectionlessideal}
\end{lemma}

\begin{lemma}  Let $G$ be a row-finite graph with no
sinks that satisfies Condition~(L), and let $\tau:C^*(G)
\rightarrow \Q$ be an essential extension of $C^*(G)$.  If $\{
s_e, p_v \}$ is the canonical Cuntz-Krieger $G$-family, then
there exists a degenerate essential extension $t:C^*(G)
\rightarrow \Q$ such that $t(s_es_e^*)=\tau(s_es_e^*)$ for all
$e \in G^1$.
\label{texists}
\end{lemma}

\begin{proof}  Since $\tau$ is essential, $\{ \tau(s_es_e^*)
\}_{e \in G^1}$ is a countable set of mutually orthogonal
nonzero projections and we may use Lemma
\ref{manyprojectionlift} to lift them to a collection $\{ R_e
\}_{e \in G^1}$ of mutually orthogonal nonzero projections in
$\B$.  Now each $\Hi_e:= \text{im }R_e$ is
infinite-dimensional, and for each $v \in G^0$ we define $\Hi_v
= \bigoplus_{ \{ s(e)=v \} } \Hi_e$.  Then each $\Hi_v$ is
infinite-dimensional and for each $e \in G^1$ we can let $T_e$
be a partial isometry with initial space $\Hi_{r(e)}$ and
final space $\Hi_e$.  Also for each $v \in G^0$ we shall let
$Q_v$ be the projection onto $\Hi_v$.  Then $\{ T_e, Q_v
\}$ is a Cuntz-Krieger $G$-family.  By the universal property
of $C^*(G)$ there exists a homomorphism $\tilde{t}:C^*(G)
\rightarrow \B$ such that $\tilde{t}(p_v)=Q_v$ and
$\tilde{t}(s_e)=T_e$. Let $t := \pi \circ \tilde{t}$.  Then $t$
is a degenerate extension and $t(s_es_e^*) =
\pi(\tilde{t}(s_es_e^*)) = \pi (T_eT_e^*) =
\pi(R_e) = \tau(s_es_e^*)$. Furthermore, for all $v
\in G^0$ we have that $$ t(p_v) = \sum_{s(e)=v}t(s_es_e^*) =
\sum_{s(e)=v}\tau(s_es_e^*) = \tau(p_v) \neq 0 $$
so $t$ is essential. \end{proof}

\begin{remark}  Suppose that $G$ is a graph with no sinks,
$\tau$ is an extension of $C^*(G)$, and $t$ is another
extension for which $t(s_es_e^*)=\tau(s_es_e^*)$.  Then $t$ will
also have the property that $t(p_v)=t(\sum s_es_e^*) =
\sum t(s_es_e^*) = \sum \tau(s_es_e^*) = \tau (\sum s_es_e^*) =
\tau(p_v)$ for any $v \in G^0$.  
\label{tagree}
\end{remark}
 
\begin{definition}
Let $\tau:C^*(G) \rightarrow \Q$ be an
essential extension of $C^*(G)$, and for each $e \in G^1$ define
$E_e := \tau(s_es_e^*)$.  If $t:C^*(G) \rightarrow \Q$ is
another essential extension of $C^*(G)$ with the property that
$t(s_es_e^*) = E_e$, then we define a vector $d_{\tau,t} \in
\prod_{G^1} \mathbb{Z}$ by $$d_{\tau,t} (e) =
-\mathrm{ind}_{E_e} \tau(s_e) t(s_e^*).$$  Note that this is
well defined since $E_e \tau(s_e) t (s_e^*) E_e = \tau (s_e) t
(s_e^*)$ and by Remark~\ref{tagree} we have that
$ \tau(s_e) t (s_e^*)  \tau(s_e^*) t (s_e)
 = \tau(s_e) \tau(s_e^*s_e) \tau(s_e^*) = E_e$ so
$\tau(s_e) t (s_e^*)$ is invertible in $E_e \Q E_e$.
\end{definition}

\begin{remark}  If $E \in \Q$ is a projection and $E' \in \B$
is a lift of $E$ to a projection in $\B$, then one can see that
$\Q(E'(\Hi)) \cong E \Q E$ via the obvious correspondence.  In
the rest of this paper we shall often identify $\Q(E'(\Hi))$
with $ E \Q E$.
\label{identifycalkinprojection}
\end{remark}

\noindent The proof of the following lemma is straightforward.

\begin{lemma}
Let $E \in \Q$ be a projection and $X \in \Q$, and suppose
that $EXE$ is invertible in $E \Q E$.  If $V \in \Q$ is a
partial isometry with initial projection $V^*V=E$ and final
projection $VV^* = F$, then $\ind_EX = \ind_F VXV^*$.
\label{usefulindexlemma}
\end{lemma}

\begin{proposition}  Let $G$ be a row-finite graph with no
sinks that satisfies Condition~(L).  Also let
$\tau$ be an essential extension of $C^*(G)$ and
$E_e:=\tau(s_es_e^*)$ for $e \in G^1$.  If $t$ and $t'$ are
essential extensions of $C^*(G)$ that are
CK-equivalent and satisfy $t(s_es_e^*)=t'(s_es_e^*) = E_e$,
then $d_{\tau,t}- d_{\tau,t'} \in \im (B_G-I)$.
\label{sameclassincoker}
\end{proposition}

\begin{proof}  Since $t$ and $t'$
are CK-equivalent, there exists a partial isometry $U
\in \Q$ such that $t = \mathrm{Ad}(U) \circ t'$ and $t' =
\mathrm{Ad}(U^*) \circ t$.  Now notice that $U$ commutes
with $E_e$.  Thus for any $e \in G^1$ we have
$\tau(s_es_e^*) = \sum_{s(f)=r(e)}
\tau(s_fs_f^*)=\sum_{s(f)=r(e)} t(s_fs_f^*) = t(s_e^*s_e)$ and
\begin{align*}
d_{\tau,t}(e) - d_{\tau,t'}(e) &= - \text{ind}_{E_e} \tau(s_e)
t(s_e^*) + \text{ind}_{E_e} \tau(s_e)t'(s_e^*) \\
&= \text{ind}_{E_e} t(s_e) \tau(s_e^*) + \text{ind}_{E_e}
\tau(s_e)t'(s_e^*) \\
&= \text{ind}_{E_e} t(s_e)\tau(s_e^*s_e) t'(s_e^*) \quad \quad
\text{by Lemma~\ref{CKlemma2}} \\
&= \text{ind}_{E_e} t(s_e)t'(s_e^*) \\
&= -d_{t,t'}(e).
\end{align*}
Hence $d_{\tau,t} - d_{\tau,t'} = -d_{t,t'}$.  Now let $k \in
\prod_{G^1} \Z$ be the vector given by $k(f) :=
\text{ind}_{E_f} U$.  Then for any $e \in G^1$ we have 
\begin{align*}
d_{t,t'}(e) &= -\text{ind}_{E_e}t(s_e)t'(s_e^*) \\
&= -\text{ind}_{E_e} t(s_e)Ut(s_e^*)U^* \\
&= -\text{ind}_{E_e} t(s_e)Ut(s_e^*)-\text{ind}_{E_e} U^* \quad
\quad \text{by Lemma~\ref{CKlemma2}} \\
&= -\text{ind}_{t(s_e^*s_e)} U -\text{ind}_{E_e}U^* \quad \quad
\text{by Lemma~\ref{usefulindexlemma}} \\
&= -\ind_{\sum E_f \atop s(f)=r(e)} U + \text{ind}_{E_e}U \\
&= -\sum_{s(f)=r(e)} \text{ind}_{E_f}U + \text{ind}_{E_e}U
\quad \quad \text{by Lemma~\ref{CKlemma1}} \\
&= - \left( \sum_{f \in G^1} B_G(e,f) k(f) - k(e) \right)
\end{align*}
so $d_{t,t'} = -(B_G-I)k$ and $d_{\tau,t}-d_{\tau,t'}=-d_{t,t'}
\in \im (B_G-I)$.
\end{proof}

\begin{definition}
Let $G$ be a row-finite graph with no sinks that satisfies
Condition~(L).  Let $B_G$ be the edge matrix of $G$ and
$B_G-I: \prod_{G^1} \Z \rightarrow \prod_{G^1} \Z$.  If
$\tau$ is an essential extension of $C^*(G)$, then we shall
define an element $d_\tau \in \coker (B_G-I)$ by $$d_\tau :=
[d_{\tau,t }] \in \coker(B_G-I),$$ where $t$ is any degenerate
extension with the property that $t(s_es_e^*)=\tau(s_es_e^*)$
for all $e \in G^1$.
\end{definition}

In the above definition, the existence of $t$ follows from
Lemma \ref{texists}.  In addition, since any two degenerate
essential extensions are CK-equivalent, it follows from
Proposition \ref{sameclassincoker} that the class of
$d_{\tau,t}$ in $\coker(B_G-I)$ will be independent of the
choice of $t$.  Therefore $d_\tau$ is well defined.

The proof of the following lemma is straightforward.

\begin{lemma}  Suppose that $\tau_1$ and $\tau_2$ are
extensions of a $C^*$-algebra $A$, and that $v$ is a partial
isometry in $\Q$ for which $\tau_1 = \mathrm{Ad}(v) \circ
\tau_2$ and $\tau_2 = \mathrm{Ad}(v^*) \circ \tau_1$.  Then
there exists either an isometry or coisometry $W \in \B$ such
that $\tau_1 = \Ad \pi(W) \circ \tau_2$ and $\tau_2 =
\Ad \pi(W^*) \circ \tau_1$.
\label{isometryreplace}
\end{lemma}

\begin{corollary}
Let $\tau_1$ and $\tau_2$ be essential extensions of a
$C^*$-algebra $A$.  Then $\tau_1$ and $\tau_2$ are CK-equivalent
if and only if there exists either an isometry or coisometry
$W$ in $\B$ such that $\tau_1 = \Ad \pi(W) \circ \tau_2$ and
$\tau_2 = \Ad \pi(W^*) \circ \tau_1$.
\label{CKisometryreplace}
\end{corollary}

\begin{lemma}  Let $G$ be a row-finite graph with no sinks that
satisfies Condition~(L).  Suppose that $\tau_1$ and
$\tau_2$ are two essential extensions of $C^*(G)$ that are
equal in $\ext(C^*(G))$.  Then $d_{\tau_1}$ and $d_{\tau_2}$
are equal in $\coker (B_G-I)$.
\label{tauindependentofCKclass}
\end{lemma}
\begin{proof}  Since $\tau_1$ and $\tau_2$ are equal in
$\ext(C^*(G))$ it follows that they are CK-equivalent.  By
interchanging
$\tau_1$ and $\tau_2$ if necessary, we may use Corollary
\ref{CKisometryreplace} to choose an isometry
$W$ in $\B$ for which $\tau_1
= \Ad \pi(W) \circ \tau_2$ and $\tau_2 =
\Ad \pi(W^*) \circ \tau_1$.  For each $e \in G^1$
define $E_e := \tau_1(s_es_e^*)$ and $F_e :=
\tau_2(s_es_e^*)$.  By Lemma \ref{texists} there exists a
degenerate essential extension $t_2 = \pi \circ \tilde{t}_2$
with the property that
$t_2(s_es_e^*) = \tau_2(s_es_e^*) = F_e$ for all $e \in G^1$. 
Then $\tilde{t}_1 := W \tilde{t}_2 W^*$ will be a
representation of $C^*(G)$ ($\tilde{t}_1$ is multiplicative
since $W$ is an isometry), and thus $t_1 := \pi \circ
\tilde{t}_1$ will be a degenerate extension with the property
that $t_1(s_es_e^*) = \tau_1(s_es_e^*)$. Now since $\tau_1$ is
essential we have that 
$$ t_1(p_v) = \sum_{s(e)=v}t_1(s_es_e^*) =
\sum_{s(e)=v}\tau_1(s_es_e^*) = \tau_1(p_v) \neq 0.$$
Therefore $p_v \notin \ker t_1$ for all $v \in G^0$ and it
follows from Lemma \ref{projectionlessideal} that
$\ker t_1 = \{ 0 \}$, and thus $t_1$ is essential.

Now recall that $E_e := \tau_1(s_es_e^*)$ and $F_e :=
\tau_2(s_es_e^*)$.  Since $W$ is an isometry, we see that
$\pi(W)F_e$ is a partial isometry with source projection $F_e$
and range projection $E_e$.  Therefore by Lemma
\ref{usefulindexlemma} it follows that 
\begin{align}
\ind_{F_e} \tau_2(s_e)t_2(s_e^*) = & \ \ind_{E_e} \pi (W) F_e
\tau_2(s_e) t_2(s_e^*) F_e \pi(W^*) \notag \\
= & \ \ind_{E_e} \pi(W) \tau_2(s_e) t_2(s_e^*) \pi(W^*) \notag
\\
= & \ \ind_{E_e} \pi(W) \tau_2(s_e) \pi(W^*) \pi(W) t_2(s_e^*)
\pi(W^*) \notag \\
= & \ \ind_{E_e}  \tau_1(s_e) t_1(s_e^*) \notag 
\end{align}
and $d_{\tau_2}$ equals $d_{\tau_1}$ in $\coker(B_G-I)$.
\end{proof}

\begin{definition}
If $G$ is a row-finite graph with no sinks that satisfies
Condition~(L), we define the \emph{Cuntz-Krieger map} to be the
map $d:\ext(C^*(G)) \rightarrow \coker(B_G-I)$ defined by
$\tau \mapsto d_\tau$.
\end{definition}

The previous lemma shows that the Cuntz-Krieger map $d$ is
well defined, and the next lemma shows that it is a
homomorphism.

\begin{lemma}  Suppose that $G$ is a row-finite graph with no
sinks that satisfies Condition~(L).  Then the Cuntz-Krieger
map is additive.
\label{disadditive}
\end{lemma}

\begin{proof} Let $\tau_1 $ and $\tau_2$ be elements
of $\ext (C^*(G))$ and choose the representatives $\tau_1$ and
$\tau_2$ such that $\tau_1 \perp
\tau_2$.  Let $t_1$ and $t_2$ be degenerate essential extensions
such that $t_1(s_es_e^*)=\tau_1(s_es_e^*)$ and
$t_2(s_es_e^*)=\tau_2(s_es_e^*)$.  

If we let $t = t_1 \boxplus t_2$, then it is straightforward to
see that $d_{\tau_1 \boxplus \tau_2, t} = d_{\tau_1,t_1} +
d_{\tau_2,t_2}$.  Also since $\tau_1 \boxplus \tau_2$ is weakly
stably equivalent to $\tau_1 + \tau_2$,
Lemma~\ref{tauindependentofCKclass} implies that  we have
$d_{\tau_1 \boxplus \tau_2} = d_{\tau_1 + \tau_2}$ in
$\coker(B_G-I)$.  Putting this all together gives $d_{\tau_1 +
\tau_2} = d_{\tau_1 \boxplus \tau_2} = [d_{\tau_1 \boxplus
\tau_2,t}] = [d_{\tau_1,t_1} + d_{\tau_2,t_2}] =
[d_{\tau_1,t_1}] + [d_{\tau_2,t_2}] = d_{\tau_1} + d_{\tau_2}$
in $\coker (B_G-I)$.  Thus $d$ is additive.  \end{proof}

\noindent We mention the following lemma whose proof is
straightforward.

\begin{lemma}
Let $E \in \Q$ be a projection, and suppose that $T$ is a
unitary in $E \Q E$ with $\ind_E T =0$.  If $E' \in
\B$ is a projection such that $\pi(E')=E$, then there is a
unitary $U
\in \B (E' \Hi )$ such that $\pi (U) = T$.
\label{injectivetool}
\end{lemma}

\begin{proposition}  Let $G$ be a row-finite graph with no
sinks that satisfies Condition~(L).  Then the Cuntz-Krieger
map $d : \mathrm{Ext}(C^*(G)) \rightarrow \coker (B_G-I)$
defined by $ \tau  \mapsto  d_\tau $ is injective.
\label{disinjective}
\end{proposition}

\begin{proof} Let $\tau$ be an essential extension of $C^*(G)$
and suppose that $d_\tau$ equals $0$ in $\coker (B_G-I)$.  Use
Lemma \ref{texists} to choose a degenerate essential extension
$t : = \pi \circ \tilde{t}$ of $C^*(G)$ such that
$t(s_es_e^*) = E_e := \tau(s_es_e^*)$ for all $e \in
G^1$.  Also let $E_e' := \tilde{t}(s_es_e^*)$.  

By hypothesis, there exists $k \in \prod_{G^1} \Z$ such that
$d_{\tau,t} = (B_G-I)k$.  Since $\tau$ is essential, for all $e
\in G^1$ we must have that $\pi(E_e')= E_e = \tau(s_es_e^*)
\neq 0$.  Since $E_e'$ is a projection, this implies that
$\text{dim}(\text{im }(E_e')) =
\infty$.  Therefore for each $e \in G^1$ we may choose
isometries or coisometries $V_e$ in $\B(E_e'(\Hi))$ such that
$\text{ind}_{E_e}V_e=-k(e)$.  Extend each $V_e$ to all of $\Hi$
by defining it to be zero on $(E_e'(\Hi))^\perp$.  Let $U :=
\sum_{e \in G^1} V_e$.  It follows that this sum converges in
the strong operator topology.  Notice that for all $e, f \in
G^1$ we have $$V_f \tilde{t}(s_es_e^*) = V_f E_f' E_e' =
\begin{cases} V_f &
\text{if $e=f$} \\ 0 & \text{otherwise.} \end{cases}$$ Since
$U$ commutes with $E_e'$ for all $e \in G^1$, we see that
$\pi(U) \tau(s_e) \pi (U^*) t(s_e^*)$ is a unitary in $E_e \Q
E_e$.  Hence we may consider
$\ind_{E_e} \pi(U) \tau(s_e) \pi (U^*) t(s_e^*).$  Using the
above identity we see that for each $e \in G^1$ we have
\begin{align}
\text{ind}_{E_e}\pi(U)\tau(s_e)\pi(U^*)t(s_e^*) 
= & \ \text{ind}_{E_e}\pi(U)\tau(s_es_e^*)\tau(s_e)\pi(U^*)
t(s_e^*) \notag \\ 
= & \ \text{ind}_{E_e}\pi(V_e) \tau(s_e)t(s_e^*) \Big(
t(s_e)\pi(U^*)t(s_e^*) \Big). \label{anothersumproductthing}
\end{align}
Now since $t(s_e)$ is a partial isometry with source
projection $$t(s_e^*s_e)= \sum_{s(f)=r(e)}t(s_fs_f^*) =
\sum_{s(f)=r(e)} E_f$$ and range projection $t(s_es_e^*)=E_e$,
we may use Lemma \ref{usefulindexlemma} to conclude that 
$$\ind_{{\sum E_f} \atop {s(f)=r(e)}} \pi(U^*) = \ind_{E_e}
t(s_e) \pi(U^*) t(s_e^*).$$  This combined with Lemma
\ref{CKlemma1} implies that
\begin{align}
\ind_{E_e} t(s_e) \pi(U^*) t(s_e^*) = & \ \sum_{s(f)=r(e)}
\ind_{E_f} \pi(U^*) \notag \\
= & \ \sum_{s(f)=r(e)} \ind_{E_f} \pi(V_f^*) 
\notag \\
= & \ \sum_{s(f)=r(e)} k(f) \notag \\
= & \ \sum_{f \in G^1} B_G(e,f) k(f).
\label{onceagainasumproductthing}
\end{align}
Combining (\ref{anothersumproductthing}) and
(\ref{onceagainasumproductthing}) with Lemma
\ref{CKlemma2} gives  $$\ind_{E_e} \pi(U) \tau(s_e) \pi (U^*)
t(s_e^*) = \left( \sum_{f \in G^1}B_G(e,f) k(f) - k(e)
\right) - d_\tau(e) = 0.$$
Thus by Lemma \ref{injectivetool} there exists an operator $X_e
\in \B$ such that the restriction of $X_e$ to $E_e'(\Hi)$ is a
unitary operator and $\pi
(X_e)=\pi(U)\tau(s_e)\pi(U^*)t(s_e^*)$.  Let
$T_e := X_e \tilde{t}(s_e)$.  Then $T_e$ is a partial isometry
that satisfies $T_eT_e^* = E_e'$ and
$T_e^*T_e=\tilde{t}(s_e^*)X_e^*X_e\tilde{t}(s_e) =
\tilde{t}(s_e^*s_e) = \tilde{t}(p_{r(e)})$.  One can then
check that $\{
\tilde{t}(p_v), T_e \}$ is a Cuntz-Krieger $G$-family in $\B$. 
Thus by the universal property of $C^*(G)$ there exists a
homomorphism $\tilde{\rho} :C^*(G) \rightarrow \B$ such that
$\tilde{\rho}(p_v) = \tilde{t}(p_v)$ and
$\tilde{\rho}(s_e)=T_e$.  Let $\rho := \pi \circ
\tilde{\rho}$.  Then $\rho$ is a degenerate extension of
$C^*(G)$.  Furthermore, since $\rho (p_v) = t (p_v) \neq 0$ we
see that $p_v \notin \ker \rho$ for all $v \in G^0$.  Since $G$
satisfies Condition~(L), it follows from Lemma
\ref{projectionlessideal} that $\ker \rho = \{ 0 \}$ and $\rho$
is a degenerate essential extension.
In addition, we see that for each $e \in G^1$
\begin{align}
\rho(s_e) = & \ \pi (T_e) \notag \\
= & \ \pi (X_e \tilde{t} (s_e) ) \notag \\
= & \ \pi (U) \tau(s_e) \pi (U^*) t(s_e^*) t(s_e) \notag \\
= & \ \pi (U) \tau (s_e) \pi (U^*). \notag 
\end{align} 
Thus $\rho(s_e) = \pi(U) \tau(s_e) \pi(U^*)$ for all $e \in
G^1$, and since the $s_e$'s generate $C^*(G)$, it follows that
$\rho(a) = \pi(U) \tau(a) \pi(U^*)$ for all $a \in C^*(G)$ and
hence $\rho = \text{Ad}(\pi(U)) \circ \tau$.

In addition, since the $V_e$'s are either isometries or
coisometries on $E_e'(\Hi)$ with finite Fredholm index, it
follows that
$\pi(V_e^*V_e)=\pi(V_eV_e^*)=\pi(E_e')$. 
Therefore, for any $e \in G^1$ we have that
\begin{align}
\pi(U^*U)\tau(s_e) = & \ \pi \left( U^* \sum_{f \in G^1} V_f
\tilde{t}(s_es_e^*) \right) \tau(s_e) \notag \\
= & \ \pi \left( U^* V_e E_e' \right) \tau(s_e) \notag \\
= & \ \pi \left( \sum_{f \in G^1} V_f^* E_e' V_e \right)
\tau(s_e) \notag \\
= & \ \pi ( V_e^* V_e ) \tau(s_e) \notag \\
= & \ \pi ( E_e') \tau(s_e) \notag \\
= & \ \tau(s_es_e^*) \tau(s_e) \notag \\
= & \ \tau(s_e). \notag
\end{align}
Again, since the $s_e$'s generate $C^*(G)$, it follows that
$\pi(U^*U) \tau(a) = \tau(a)$ for all $a \in C^*(G)$. 
Similarly, $ \tau(a) \pi(U^*U) = \tau(a)$ for all $a \in
C^*(G)$.  Thus $\pi(U^*) \rho (a) \pi(U) = \pi(U^*U) \tau(a)
\pi(U^*U) = \tau(a)$ for all $a \in C^*(G)$ and $\tau =
\text{Ad}(\pi(U)^*) \circ \rho$.

Now because the $V_e$'s are all isometries or coisometries on
orthogonal spaces, it follows that $U$, and hence $\pi(U)$, is
a partial isometry.  Therefore, $\tau = \rho$ in $\ext (C^*(G))$
and since $\rho$ is a degenerate essential extension it follows
that $\tau = 0$ in $\ext (C^*(G))$.  This implies that $d$ is
injective.
\end{proof}

\section{The $\W$ Map}
\label{Wmapsection}

In the previous section we showed that if $G$ is a row-finite
graph with no sinks that satisfies Condition~(L), then the
Cuntz-Krieger map $d : \ext(C^*(G)) \rightarrow \coker(B_G-I)$
is a monomorphism.  It turns out that $d$ is also surjective;
that is, it is an isomorphism.  In this section we shall prove
this fact, but we shall do it in an indirect way.  We show
that $\coker (B_G-I)$ is isomorphic to $\coker (A_G-I)$ and
then compose $d$ with this isomorphism to get a map from
$\ext(C^*(G))$ into
$\coker (A_G-I)$.  We call this composition the $\W$ map
and we shall show that it, and consequently also $d$, is
surjective.  For the rest of this paper we will be mostly
concerned with the $\W$ map and how it relates to 1-sink
extensions defined in \cite{RTW}.

\begin{definition}
Let $G$ be a graph.  The \emph{source matrix} of $G$ is the $G^0
\times G^1$ matrix given by $$S_G(v,e) = \begin{cases} 1 &
\text{if $s(e)=v$} \\ 0 & \text{otherwise} \end{cases}$$
and the \emph{range matrix} of $G$ is the $G^1
\times G^0$ matrix given by $$R_G(e,v) = \begin{cases} 1 &
\text{if $r(e)=v$} \\ 0 & \text{otherwise.} \end{cases}$$
\end{definition}
Notice that if $G$ is a row-finite graph, then $S_G$ will have
rows that are eventually zero and left multiplication by
$S_G$ defines a map $S_G :
\prod_{G^1} \Z
\rightarrow \prod_{G^0} \Z$.  Also $R_G$
will always have rows that are eventually zero.  (In fact,
regardless of any conditions on
$G$, $R_G$ will have only one nonzero entry in each row.) 
Therefore left multiplication by $R_G$ defines a map $R_G :
\prod_{G^0} \Z
\rightarrow \prod_{G^1} \Z$.  Furthermore, one can see that
$$R_GS_G = B_G \hspace{.1in} \text{ and } \hspace{.1in}
S_GR_G=A_G.$$

The following lemma is well known for finite graphs and a proof
for $S_G$ restricted to the direct sum $S_G : \bigoplus_{G^1} \Z
\rightarrow \bigoplus_{G^0} \Z$ is given in \cite[Lemma
4.2]{MRS}.  Essentially the same proof goes through if we
replace the direct sums by direct products.

\begin{lemma}
Let $G$ be a row-finite graph.  The map $S_G : \prod_{G^1} \Z
\rightarrow \prod_{G^0} \Z$ induces an
isomorphism $\overline{S_G} : \coker (B_G-I) \rightarrow \coker
(A_G-I)$.
\label{Sisomorphism}
\end{lemma}

\begin{proof} Suppose that $z \in \im (B_G-I)$.  Then
$z=(B_G-I)u$ for some $u \in \prod_{G^1} \Z$.  Then 
$$S_G z = S_G (B_G-I)u = S_G (R_GS_G - I) u = (S_GR_G - I)S_G
u = (A_G-I) S_G u$$
and $S_G$ does in fact map $\im(B_G-I)$ into $\im (A_G-I)$. 
Thus $S_G$ induces a homomorphism $\overline{S_G}$ of $\coker
(B_G-I)$ into $\coker (A_G-I)$.
In the same way, $R_G$ induces a homomorphism $\overline{R_G}$
from $\coker (A_G-I)$ into $\coker(B_G-I)$, which we claim is
an inverse for $\overline{S_G}$.  We see that
\begin{align}
\overline{R_G} \circ \overline{S_G} (u + \im (B_G-I)) = &
\ R_GS_Gu + \im (B_G-I) \notag \\
= & \ u + (B_Gu-u) + \im (B_G-I) \notag \\
= & \ u + \im (B_G-I) \notag
\end{align}
and similarly $\overline{S_G} \circ \overline{R_G}$ is the
identity on $\coker (A_G-I)$. \end{proof}

\begin{definition}
Let $G$ be a row-finite graph with no
sinks that satisfies Condition~(L), and let $d: \ext(C^*(G))
\rightarrow \coker (B_G-I)$ be the Cuntz-Krieger map.  The
\emph{$\W$ map} is the homomorphism $\wmap : \ext(C^*(G))
\rightarrow
\coker(A_G-I)$ given by $\wmap := \overline{S_G} \circ d$. 
Given an extension
$\tau$ of $C^*(G)$, we shall refer to the class $\wmap (\tau)$
in $\coker (A_G-I)$ as the \emph{$\W$ class} of $\tau$.
\end{definition}

\begin{lemma}
Let $G$ be a row-finite graph with no sinks that satisfies
Condition~(L).  Then the $\W$ map associated to $G$ is a
monomorphism.
\end{lemma}
\begin{proof}  Since $\wmap = \overline{S_G} \circ d$, and
$\overline{S_G}$ is an isomorphism by Lemma
\ref{Sisomorphism}, the result follows from
Proposition~\ref{disinjective}. \end{proof}

We shall eventually show that the $\W$ map is also surjective;
that is, it is an isomorphism.  In order to do this
we consider 1-sink extensions, which were
introduced in \cite{RTW}, and describe a way to associate
elements of $\ext(C^*(G))$ to them.

\begin{definition} \cite[Definition 1.1]{RTW} Let $G$ be a
row-finite graph.  A \emph{$1$-sink extension} of $G$ is a
row-finite graph $E$ that contains $G$ as a subgraph and
satisfies:
\begin{enumerate}
\item
$H:=E^0\setminus G^0$ is finite, contains no sources,
and contains exactly $1$ sink $v_0$.
\item
There are no loops in $E$ whose vertices lie in
$H$.
\item
If $e \in E^1 \setminus G^1$, then $r(e) \in H$.
\item
If $w$ is a sink in $G$, then $w$ is a sink in $E$.
\end{enumerate}  We will write $(E,v_0)$ for the 1-sink
extension, where $v_0$ denotes the sink outside $G$.
\end{definition}

If $(E,v_0)$ is a 1-sink extension of $G$, then we may let
$\pi_E :C^*(E) \rightarrow C^*(G)$ be the surjection described
in \cite[Corollary 1.3]{RTW}.  Then $\ker \pi_E = I_{v_0}$
where $I_{v_0}$ is the ideal in $C^*(E)$ generated by the
projection $p_{v_0}$.  Thus we have a
short exact sequence 
\begin{equation}
\xymatrix{ 0 \ar[r] & I_{v_0} \ar[r]^i & C^*(E) \ar[r]^{\pi_E}
& C^*(G) \ar[r] & 0.}
\notag
\end{equation}
We call $E$ an \emph{essential} 1-sink extension of
$G$ when $G^0 \geq v_0$.  Note that $I_{v_0}$ is an
essential ideal of $C^*(E)$ if and only if $E$ is an
essential 1-sink extension of $G$ \cite[Lemma 2.2]{RTW}. 

\begin{lemma}
If $G$ is a row-finite graph and $(E,v_0)$ is an
essential 1-sink extension of $G$, then $I_{v_0} \cong \K$.
\label{idealisthecompacts}
\end{lemma}

\begin{proof}  Let $E^*(v_0)$ be the set of all paths in
$E$ whose range is $v_0$.  Since
$E$ is an essential 1-sink extension of
$G$, it follows that $G^0 \geq v_0$.  Thus
for every $w \in G^0$ there exists a path from $w$ to $v_0$.  If
$G^0$ is infinite, this implies that $E^*(v_0)$ is
also infinite.  If $G^0$ is finite, then because $G^0 \geq v_0$
it follows that $G$ is a finite graph with no sinks, and hence
contains a loop.  If $w$ is any
vertex on this loop, then there is a path from $w$ to $v_0$ and
hence $E^*(v_0)$ is infinite.  Now because $E^*(v_0)$ is
infinite it follows from \cite[Corollary 2.2]{KPR} that
$I_{v_0} \cong \K ( \ell^2(E^*(v_0))) \cong \K$. \end{proof}

\begin{definition}
Let $G$ be a row-finite graph and let $(E,v_0)$ be an essential
1-sink extension of
$G$.  \emph{The extension associated to
$E$} is (the strong equivalence class of) the
Busby invariant of any extension
\begin{equation}
\xymatrix{ 0 \ar[r] & \K \ar[r]^{i_E} & C^*(E) \ar[r]^{\pi_E}
& C^*(G) \ar[r] & 0}
\notag
\end{equation}
where $i_E$ is any isomorphism from $\K$ onto $I_{v_0}$. 
As with other extensions
we shall not distinguish between an extension and its Busby
invariant.
\end{definition}

\begin{remark}
The above extension is well-defined up to strong
equivalence.  If different choices of $i_E$ are made then it
follows from a quick diagram chase that the two associated
extensions will be strongly equivalent (see problem 3E(c) of
\cite{WO} for more details).  Also recall that since $p_{v_0}$
is a minimal projection in $I_{v_0}$ \cite[Corollary 2.2]{KPR},
it follows that $i_E^{-1} (p_{v_0})$ will always be a rank 1
projection in $\K$.
\end{remark}

Let $(E,v_0)$ be a 1-sink extension of $G$.  Then for
$w \in E^0$ we denote by $Z(w,v_0)$ the set of paths $\alpha$
from $w$ to $v_0$ with the property that $\alpha_i \in E^1
\backslash G^1$ for $1 \leq i \leq | \alpha |$.  The
\emph{$\W$ vector} of $E$ is the element $\wvec_E \in
\prod_{G^0} \N$ given by
$$\wvec_E(w) := \# Z(w,v_0).$$  An edge $e \in E^1$ with
$s(e) \in G^0$ and $r(e) \notin G^0$ is called a \emph{boundary
edge}, and the sources of these edges are called \emph{boundary
vertices}.  

\begin{lemma}
\label{usefullemmaonrankforrepresentations}
Let $G$ be a row-finite graph and let $(E,v_0)$ be a
1-sink extension of $G$.  If
$\{s_e, p_v \}$ is the canonical Cuntz-Krieger
$E$-family in $C^*(E)$ and $\sigma : C^*(E)
\rightarrow \B$ is a representation with the property that
$\sigma ( p_{v_0})$ is a rank 1 projection, then 
$$ \rank \sigma (s_e) = \# Z(r(e),v_0) \hspace{.3in} \text{for
all $e \in E^1
\backslash G^1$.} $$
\end{lemma}
\begin{proof} For $e \in E^1 \backslash G^1$ let $k_e := \max \{
| \alpha | : \alpha \in Z((r(e),v_0) \}$.  Since $E$ is
a 1-sink extension of $G$ we know that $k_e$ is finite.  We
shall prove the claim by induction on $k_e$.  If $k_e=0$, then
$r(e)=v_0$ and $\rank \sigma(s_e) = \rank \sigma(s_e^*s_e) =
\rank \sigma (p_{v_0}) = 1$.

Assume that the claim holds for all $f \in E^1 \backslash G^1$
with $k_f \leq m$.  Then let $e \in E^1 \backslash G^1$ with
$k_e = m+1$.  Since $E$ is a 1-sink extension of $G$ there are
no loops based at $r(e)$.  Thus $k_f \leq m$ for all $f \in E^1
\backslash G^1$ with $s(f) = r(e)$.  By the induction hypothesis
$\rank \sigma (s_f) = \# Z(r(e),v_0)$ for all $f$ with $s(f) =
r(e)$.  Since the projections $s_fs_f^*$ are mutually
orthogonal we have 
\begin{align*}
\rank \sigma(s_e) & = \rank \sigma(s_e^*s_e) =
\rank \left( \sum_{ s(f) = r(e)} \sigma(s_fs_f^*) \right) =
\sum_{s(f)=r(e)}
\rank \sigma (s_fs_f^*) \\ & = \sum_{s(f)=r(e)} \# Z((r(f),v_0)
= \# Z(r(e),v_0).
\end{align*}
\end{proof}

\begin{lemma}
Let $G$ be a row-finite graph with no sinks that
satisfies Condition~(L), and let $d: \ext (C^*(G)) \rightarrow
\coker (B_G-I)$ be the Cuntz-Krieger map.  If $(E,v_0)$ is an
essential 1-sink extension of $G$ and $\tau$ is the Busby
invariant of the extension associated to $E$, then $$d (\tau) =
[x]$$ where $[x]$ is the class in $\coker(B_G-I)$ of the
vector $x \in \prod_{G^1} \Z$ given by $x(e) := \wvec_E(r(e))$
for all $e \in G^1$, and $\wvec_E$ is the $\W$ vector of $E$.
\label{dw}
\end{lemma} 

\begin{proof}  Let $\{s_e, p_v \}$ be the
canonical Cuntz-Krieger $G$-family in $C^*(G)$, and let $\{
t_e, q_v \}$ be the canonical Cuntz-Krieger $E$-family in
$C^*(E)$.  Choose an
isomorphism $i_E : \K \rightarrow I_{v_0}$, and let $\sigma$ and
$\tau$ be the homomorphisms that make the diagram
\begin{equation}
  \notag
  \xymatrix{0 \ar[r] & \K \ar@{=}[d] \ar[r]^{i_E} &
C^*(E) \ar[d]^{\sigma} \ar[r]^{\pi_E} & C^*(G) \ar[r]
\ar[d]^{\tau} & 0 \\  0 \ar[r] & \K \ar[r]^i &
\B \ar[r]^{\pi} & \Q \ar[r] & 0}
\end{equation} 
commute.  Then $\tau$ is the Busby invariant of the extension
associated to $E$, and since $E$ is an essential
1-sink extension, it follows that $\sigma$ and $\tau$ are
injective.  For all $v \in E^0$ and $e \in E^1$ define 
$$H_v := \im \sigma(q_v) \text{ \hspace{.2in} and
\hspace{.2in}} H_e := \im \sigma (t_et_e^*).$$  Note that
$s(e)=v$ implies that $H_e \subseteq H_v$.  Also since
$i_E^{-1}(q_{v_0})$ is a rank 1 projection, and since
the above diagram commutes, it follows that $\sigma(q_{v_0})$
is a rank 1 projection.  Thus $H_{v_0}$ is 1-dimensional. 
Furthermore, by Lemma
\ref{usefullemmaonrankforrepresentations} we see that
$\text{dim} (H_v) = \# Z(v,v_0)$ and $\text{dim} (H_e) = \#
Z(r(e),v_0)$ for all $v \in E^0 \backslash G^0$ and $e \in E^1
\backslash G^1$.  In addition, since $t_et_e^* \leq q_{s(e)}$
for any $e \in E^1 \backslash G^1$ and because the $q_v$'s are
mutually orthogonal projections, it follows that the $H_e$'s
are mutually orthogonal subspaces for all $e \in E^1 \backslash
G^1$. 

For all $v \in G^0$ define
$$ V_v := H_v \ominus \big( \bigoplus_{e \text{ is a boundary}
\atop \text{ edge and } s(e) =v} H_e \big).$$
Then for every $v \in G^0$, we have $\pi(\sigma(q_v)) =
\tau(\pi_E(q_v))=\tau(p_v) \neq 0$ since
$\tau$ is injective.  Therefore, the rank of
$\sigma(q_v)$ is infinite and hence $\text{dim}(H_v) = \infty$
and $\text{dim}(V_v) = \infty$.  Now for each $v \in G^0$ and $e
\in G^1$ let $P_v$ be the projection onto $V_v$ and $S_e$ be a
partial isometry with initial space $V_{r(e)}$ and final space
$H_e$.  One can then check that $\{ S_e, P_v \}$ is a
Cuntz-Krieger $G$-family in $\B$.  Therefore, by the universal
property of $C^*(G)$ there exists a homomorphism $\tilde{t} :
C^*(G) \rightarrow \B$ with the property that $\tilde{t} (s_e)
= S_e$ and $\tilde{t} (p_v) = P_v$.  Define $t := \pi \circ
\tilde{t}$.

Then for all $v \in
G^0$ we have that $$t(p_v) = \pi (\tilde{t} (p_v)) = \pi (P_v)
\neq 0.$$  Thus $p_v \notin \ker t$ for all $v \in G^0$.  By
Lemma \ref{projectionlessideal} it follows that $\ker t = \{ 0
\}$ and $t$ is an essential extension of $C^*(G)$.
Now since $S_eS_e^*$ is a projection onto a subspace of $\im
\sigma(t_et_e^*)$ with finite codimension, it follows that
$\pi (S_eS_e^*) = \pi (\sigma(t_et_e^*))$.  Thus $t$ has the
property that for all $e \in G^1$ 
$$t(s_es_e^*) = \pi(\tilde{t}(s_es_e^*)) = \pi (S_eS_e^*) =
\pi ( \sigma ( t_e t_e^*)) = \tau (\pi_E(t_et_e^*)) =
\tau(s_es_e^*).$$  By the definition of the Cuntz-Krieger map
$d$ it follows that the image of the extension associated to $E$
will be the class of the vector $d_\tau$ in $\coker (B_G-I)$,
where
$d_\tau(e) = - \text{ind}_{\tau(s_es_e^*)}\tau(s_e)t(s_e^*)$.
Now $\text{ind}_{\tau(s_es_e^*)}\tau(s_e)t(s_e^*)$ is equal to
the Fredholm index of
$\sigma(t_et_e^*)\sigma(t_e)S_e^*\sigma(t_et_e^*)=
\sigma(t_e)S_e^* \ $ in $\im (\sigma(t_et_e^*))=H_e$.
Since $S_e$ is a partial isometry with initial space
$V_{r(e)} \subseteq H_{r(e)}$ and final space $H_e$,
and since $\sigma(t_e)$ is a partial isometry with initial
space $H_{r(e)}$ it follows that $\ker \sigma(t_e)S_e^* = \{ 0
\}$ in $H_e$.  Furthermore, $\sigma(t_e^*)$ is a partial
isometry with initial space $H_e$ and final space $$H_{r(e)} =
V_{r(e)} \oplus \big( \bigoplus_{f \text{ is a boundary}
\atop \text{ edge and } s(f) =r(e)} H_f \big)$$
and $S_e$ is a partial isometry with initial space
$V_{r(e)}$.  Therefore, since $\text{dim}(H_f)= \# Z(r(f),v_0)$
for all $f \notin G^1$ we have that 
$$ \ker((\sigma(t_e)S_e)^*) = \ker(S_e \sigma(t_e^*))
= \sum_{s(f)=r(e)} Z(r(f),v_0) = \wvec_E(r(e)).$$ Thus
$d_\tau(e) = \wvec_E(r(e))$ for all $e \in G^1$.
\end{proof}

\begin{proposition}
\label{wclassequalswvec}
Let $G$ be a row-finite graph with no sinks that satisfies
Condition~(L), and suppose that $(E,v_0)$ is an essential 1-sink
extension of $G$.  If $\tau$ is the Busby invariant of the
extension associated to $E$, then the value that the $\W$ map
$\wmap: \ext(C^*(G)) \rightarrow \coker(A_G-I)$ assigns to
$\tau$ is given by the class of the $\W$ vector in
$\coker (A_G-I)$; that is, 
$$\wmap(\tau) = [ \wvec_E ].$$
\end{proposition} 
\begin{proof} From Lemma \ref{dw} we have that $d_\tau = [x]$
in $\coker(B_G-I)$, where $x \in \prod_{G^1} \Z$ is the vector
given by $x(e) :=
\wvec_E(r(e))$ for $e \in G^1$.  By the definition of $\wmap$
we have that $\wmap(\tau) := \overline{S_G}(d_\tau)$ in
$\coker(A_G-I)$.  Thus $\wmap (\tau)$ equals the class of the
vector $y \in \prod_{G^0} \Z$ given by $$y(v) = (S_G(x))(v) =
\sum_{s(e)=v} x(e) = \sum_{s(e)=v} \wvec_E(r(e)).$$  Hence for
all $v \in G^0$ we have $$y(v) - \wvec_E(v) = \sum_{s(e)=v}
\wvec_E(r(e)) - \wvec_E(v) = \sum_{w \in G^0}
A_G(v,w)\wvec_E(w) - \wvec_E(v)$$so $y - \wvec_E = (A_G-I)
\wvec_E$.  Thus $[y] = [\wvec_E]$ and $\wmap(\tau) = [\wvec_E]$
in $\coker(A_G-I)$. \end{proof}

This result gives us a method to prove that $\wmap$ is
surjective.  We need only produce essential 1-sink extensions
with the appropriate $\W$ vectors.

A 1-sink extension  $E$ of $G$ is said to be \emph{simple} if
$E^0 \backslash G^0$ consists of a single vertex.  If $G$ is a
graph with no sinks, then for any $x \in \prod_{G^0} \N$ we may
form a simple 1-sink extension of $G$ with $\W$ vector equal to
$x$ merely by defining $E^0 := G^0 \cup \{ v_0 \}$ and $E^1 :=
G^1 \cup \{ e_w^i : w \in G^0 \text{ and } 1 \leq i \leq x(w)
\}$ where each $e_w^i$ is an edge with source $w$ and range
$v_0$.  In order to show that the $\W$ map is
surjective we will not only need to produce such 1-sink
extensions, but also ensure that they are essential.

\begin{lemma}
Let $G$ be a row-finite graph with no sinks that satisfies
Condition~(L).  There exists a vector $n \in \prod_{G^0} \Z$
with the following two properties:
\begin{enumerate}
\item
$(A_G-I)n \in \prod_{G^0} \N$ 
\item for all $v \in G^0$ there exists $w \in G^0$ such that
$v \geq w$ and $((A_G-I)n)(w) \geq 1$.
\end{enumerate}
\label{vectortomakeessential}
\end{lemma}
\begin{proof}  Let $L \subseteq G^0$ be those vertices of $G$
that feed into a loop; that is, 
$$L :=  \{ v \in G^0 : \text{there exists a loop $x$ in $G$
for which $v \geq r(x_1)$} \}.$$ 
Now consider the set $M := G^0 \backslash L$.  Because
$G$ has no sinks, and because $ v \in M$ and $v \geq w$
implies that $w \in M$, it follows that $M$ cannot have a
finite number of elements.  Thus $M$ is either empty or
countably infinite.  If $M \neq \emptyset$ then list the
elements of $M$ as $M = \{w_1, w_2, \ldots \}$.
Now let $v_1^1 := w_1$.  Choose an edge $e_1^1 \in G^1$ with
the property that $s(e_1^1) = v_1^1$ and define $v_2^1 :=
r(e_1^1)$.  Continue in this fashion: given $v_k^1$ choose an
edge $e_k^1$ with $s(e_k^1)=v_k^1$ and define $v_{k+1}^1 :=
r(e_k^1)$.  Then $v_1^1, v_2^1, \ldots$ are the vertices of
an infinite path which are all elements of $M$.  Since these
vertices do not feed into a loop it follows that they are
distinct; i.e. $v_i^1 \neq v_j^1$ when $i \neq j$.

Now if every element $w \in M$ has the property that $w \geq
v_i^1$ for some $i$, then we shall stop.  If not, choose the
smallest $j \in \N$ for which $w_j \ngeq v_i^1$ for all $i \in
\N$.  Then define $v_1^2 := w_j$ and choose an edge $e_1^2$
with $s(e_1^2)=v_1^2$.  Define $v_2^2 := r(e_1^2)$.  Continue
in this fashion:  given $v_k^2$ choose an edge $e_k^2$ with
$s(e_k^2)=v_k^2$ and define $v_{k+1}^2 := r(e_k^2)$.  Then we
produce a set of distinct vertices $v_1^2, v_2^2, v_3^2,
\ldots$ that lie on the infinite path $e_1^2 e_2^2 e_3^2
\ldots$.  Moreover, since $v_1^2 \ngeq v_i^1$ for all $i$ we
must have that the $v_i^2$'s are also distinct from the
$v_i^1$'s.

Continue in this manner.  Having produced an infinite path
$e_1^k e_2^k e_3^k \ldots$ with distinct vertices $v_1^k,
v_2^k, \ldots$ we stop if every element $w \in M$ has the
property that $w \geq v_i^j$ for some $1 \leq i < \infty ,
1 \leq j \leq k$.  Otherwise, we choose the smallest $l \in
\N$ such that $w_l \ngeq v_i^j$ for all $1 \leq i < \infty, 1
\leq j \leq k$.  We define $v_1^{k+1} := w_l$.  Given
$v_j^{k+1}$ we choose an edge $e_j^{k+1}$ with
$s(e_j^{k+1})=v_j^{k+1}$.  We then define $v_{j+1}^{k+1} :=
r(e_j^{k+1})$.  Thus we produce an infinite path $e_1^{k+1}
e_2^{k+1} \ldots$ with distinct vertices $v_1^{k+1},
v_2^{k+1}, \ldots$.  Moreover, since $v_1^{k+1} \ngeq v_i^j$
for all $1 \leq i < \infty, 1 \leq j \leq k$, it follows that
the $v_i^{k+1}$'s are distinct from the $v_i^j$'s for $j \leq
k$.

By continuing this process we are able to produce the
following.  For some $n \in \N \cup \{ \infty \}$ there is a set
of distinct vertices $S \subseteq M$ given by $$S = \{ v_j^k : 1
\leq j < \infty, 1 \leq k < n \}$$
with the property that $M \geq S$, and for any $v_j^k \in S$
there exists an edge $e_j^k \in G^1$ for which
$s(e_j^k)=v_j^k$ and $r(e_j^k) = v_{j+1}^k$.

Now define $$a_v = \begin{cases} 1 & \text{if $v \in L$} \\ j
& \text{if $v=v_j^k \in S$} \\ 0 & \text{otherwise.}
\end{cases}$$and let $n := (a_v) \in \prod_{G^0} \Z$.  We
shall now show that $n$ has the appropriate properties.
We shall first show that $(A_G-I)n \in \prod_{G^0} \N$.  Let
$v \in G^0$ and consider four cases.  (Throughout the
following remember that the entries of $n$ are nonnegative
integers.)

\noindent Case 1: $A_G(v,v) \geq 1$.  Then $((A_G-I)n)(v) \geq
a_v ( A_G(v,v) - 1) \geq 0.$

\noindent Case 2: $A_G(v,v)=0, v \in L$.  Since $A_G(v,v)=0$
and $v$ feeds into a loop, there must exist an edge $e \in
G^1$ with $s(e) = v$ and $r(e) \in L$.  Thus 
$$((A_G-I)n)(v) \geq a_v (A_G(v,v)-1) + a_{r(e)} A_G(v,r(e))
\geq 1(-1)+1(1) = 0.$$

\noindent Case 3: $A_G(v,v)=0, v = v_j^k \in S$.  Then there
exists an edge $e_j^k$ with $s(e_j^k)=v_j^k$ and $r(e_j^k) =
v_{j+1}^k \neq v_j^k$.  Thus
$$((A_G-I)n)(v) \geq a_v (A_G(v,v)-1) + a_{v_{j+1}^k}
A_G(v,v_{j+1}^k) \geq j(-1)+(j+1)(1) =1.$$

\noindent Case 4: $A_G(v,v)=0, v \notin L, v \notin S$.  Then
$$((A_G-I)n)(v) \geq a_v (A_G(v,v)-1) \geq 0 \cdot
(A_G(v,v)-1) = 0.$$
Therefore $(A_G-I)n \in \prod_{G^0} \N$.  

We shall now show
that for all $v \in G^0$ there exists $w \in G^0$ such that $v
\geq w$ and $((A_G-I)n)(w) \geq 1$.  If $v \notin L$, then $v
\in M$ and $v \geq v_j^k$ for some $v_j^k \in S$.  But then
there is an edge $e_j^k$  with $s(e_j^k)=v_j^k$ and
$r(e_j^k)=v_{j+1}^k \neq v_j^k$.  Thus we have that 
\begin{align}
((A_G-I)n)(v_j^k) \geq & a_{v_j^k}
(A_G(v_j^k,v_j^k)-1) + a_{v_{j+1}^k} A_G(v_j^k,v_{j+1}^k)
\notag \\ 
\geq & (j)(0-1)+(j+1)(1) =1. \notag
\end{align}
On the other hand, if $v \in L$, then $v$ feeds into a loop. 
Since $G$ satisfies Condition~(L) this loop must have an
exit.  Therefore, there exists $w \in L$ such that $v \geq w$
and $w$ is the source of two distinct edges $e,f \in G^1$,
where one of the edges, say $e$, is the edge of a loop and
hence has the property that $r(e) \in L$.  Now consider the
following three cases.

\noindent Case 1: $r(f) \notin L$.  Then $r(f) \in M$ and
hence $r(f) \geq v_j^k$ for some $v_j^k \in S$.  But then $v
\geq v_j^k$ and $((A_G-I)n)(v_j^k) \geq 1$ as above.

\noindent Case 2:  $r(f) \in L$ and $r(e) = r(f)$.  Then 
$$((A_G-I)n)(w) \geq -a_w + a_{r(f)}A_G(w,r(f)) \geq -1 +
(1)(2) = 1.$$

\noindent Case 3:  $r(f) \in L$ and $r(e) \neq r(f)$.  Then
\begin{align}
((A_G-I)n)(w) \geq & -a_w + a_{r(e)}A_G(w,r(e)) +
a_{r(f)}A_G(w,r(f)) \notag \\
\geq & -1 + (1)(1) + (1)(1) = 1. \notag
\end{align}
\end{proof}

\begin{lemma}
Let $G$ be a row-finite graph with no sinks that satisfies
Condition~(L).  Let $x \in \prod_{G^0} \N$.  Then there exists
an essential 1-sink extension $E$ of $G$ with the property
that $[\wvec_E] = [x]$ in $\coker(A_G-I)$.
\label{Eexists}
\end{lemma}
\begin{proof}  By Lemma \ref{vectortomakeessential} we see that
there exists $n \in \prod_{G^0} \Z$ with the property that
$(A_G-I)n \in \prod_{G^0} \N$ and for all $v \in G^0$ there
exists $w \in G^0$ for which $v \geq w$ and $((A_G-I)n)(w)
\geq 1$.  Since $x+(A_G-I)n \in \prod_{G^0} \N$ we may let
$E$ be a 1-sink extension of $G$ with $\W$ vector
$\wvec_E = x + (A_G-I)n$.  Let $v_0$ be the sink of $E$.  We
shall show that $E$ is essential.  Let $v \in G^0$.  Then there
exists $w \in G^0$ for which $v \geq w$ and $((A_G-I)n) \geq
1$.  But then
$\wvec_E(w) \geq ((A_G-I)n)(w) \geq 1$ and $w$ is a boundary
vertex of $E$.  Hence $v \geq w \geq v_0$ and we have shown that
$G^0 \geq v_0$.  Thus $E$ is essential, and furthermore
$[\wvec_e] = [x+(A_G-I)n] =[x]$ in $\coker(A_G-I)$.
\end{proof}

\begin{proposition}
Let $G$ be a row-finite graph with no sinks that
satisfies Condition~(L).  The $\W$ map $\wmap : \ext(C^*(G))
\rightarrow \coker (A_G-I)$ is surjective.
\end{proposition}
\begin{proof} If $x$ is any vector in $\prod_{G^0} \N$, then by
Lemma
\ref{Eexists} there exists an essential 1-sink extensions
$E$ for which $[\wvec_E] = [x]$.  If $\tau$ is the Busby
invariant of the extension associated to $E$, then by Lemma
\ref{wclassequalswvec} we have that $\wmap (\tau) =
[\wvec_{E_1}] = [x]$.  Thus $[x] \in \im \wmap$ for all $x \in
\prod_{G^0} \N$.

Now because $C^*(G)$ is separable and nuclear, it follows from
\cite[Corollary 15.8.4]{Bla} that $\ext(C^*(G))$ is a group. 
Because $\prod_{G^0} \N$ is the positive cone of $\prod_{G^0}
\Z$, and hence generates $\prod_{G^0} \Z$, the fact that $[x]
\in \im \wmap$ for all $x \in \prod_{G^0} \N$ implies
that $\im \wmap = \coker (A_G-I)$.  \end{proof}

\begin{corollary}
Let $G$ be a row-finite graph with no sinks that
satisfies Condition~(L).  The map $d : \ext(C^*(G)) \rightarrow
\coker (B_G-I)$ is surjective.
\end{corollary}
\begin{proof}  This follows from the fact that $\wmap =
\overline{S_G} \circ d$, and $\overline{S_G}$ is
an isomorphism. \end{proof}

\begin{theorem}
\label{wisaniso}
Let $G$ be a row-finite graph with no sinks that satisfies
Condition~(L).  The $\W$ map  $\wmap : \ext(C^*(G)) \rightarrow
\coker (A_G-I)$ and the Cuntz-Krieger map $d : \ext(C^*(G))
\rightarrow \coker (B_G-I)$ are isomorphisms.  Consequently,
$$\ext(C^*(G)) \cong \coker (A_G-I) \cong \coker(B_G-I).$$
\end{theorem}

\begin{remark}
Suppose that $G$ is a row-finite graph with no sinks that
satisfies Condition~(L), and that $\tau$ is an element of $\ext
(C^*(G))$ for which $\wmap ( \tau )
\in \coker (A_G-I)$ can be written as $[x]$ for some $x \in
\prod_{G^0} \N$.  Then Lemma \ref{Eexists} shows us that there
exists an essential 1-sink extension
$E$ with the property that the extension associated to $E$ is
equal to $\tau$ in $\ext(C^*(G))$.  Thus for
every $\tau \in \ext(C^*(G))$  with the property that
$\wmap(\tau) = [x]$ for $x \in \prod_{G^0} \N$, we may choose a
representative that is the extension associated
to an essential 1-sink extension.  It is natural to
wonder if this is the case for all elements of $\ext
(C^*(G))$.  It turns out that in general it is not.  To see
this let $G$ be the following infinite graph.

$$
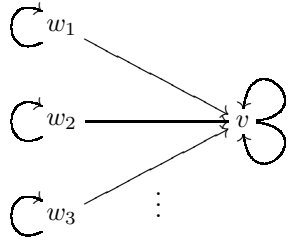
\xymatrix{  w_1 \ar@(dl,ul)[] \ar[drr] & &
\\ w_2 \ar@(dl,ul)[] \ar[rr] & & v \ar@(r,u)[]
\ar@(r,d)[] \\ 
w_3 \ar@(dl,ul)[] \ar[urr]_{\vdots} & & \\
 & & \\
}$$
\noindent Then $G$ is a row-finite graph with no sinks that
satisfies Condition (L).  However,
$$A_G-I = \left( \begin{smallmatrix} 
1 & 0 & 0 & \\
1 & 0 & 0 & \cdots \\
1 & 0 & 0 & \\
 &\vdots &   & \ddots \end{smallmatrix} \right), $$
and if we let $x : = \left( \begin{smallmatrix} -1 \\ -2 \\ -3
\\ \vdots \end{smallmatrix} \right) \in \prod_{G^0} \Z$ then
for all $n \in \prod_{G^0} \Z$ we have that
$$x + (A_G-I)n = \left( \begin{smallmatrix} -1 +n(v) \\ -2 +
n(v) \\ -3 + n(v) \\ \vdots \end{smallmatrix} \right).$$  Thus
for any $n \in \prod_{G^0} \Z$ we see that $x+(A_G-I)n$ has
negative entries.  Hence $x+(A_G-I)n$ cannot be the $\W$ vector
of a 1-sink extension for any $n \in \prod_{G^0} \Z$.  

It turns out, however, that if we add the condition that $G$ be
a finite graph then the result does hold.
\end{remark}

\begin{lemma}
Let $G$ be a finite graph with no sinks that satisfies
Condition~(L).  If $v \in G^0$, then there exists $n \in
\prod_{G^0} \N$ for which $(A_G-I)n \in \prod_{G^0} \N$ and
$((A_G-I)n)(v) \geq 1$.
\label{finitevectorexists}
\end{lemma}
\begin{proof}  If $A_G(v,v) \geq 2$ then we can let $n =
\delta_v$ and the claim holds.  Therefore, we shall suppose
that $A_G(v,v)
\leq 1$.  Since $G$ has no sinks and satisfies Condition (L),
there must exist an edge $e_1 \in G^1$ with $s(e_1)=v$ and
$r(e_1) \neq v$.  Then since $G$ has no sinks we may find an
edge $e_2 \in G^1$ with $s(e_2)=r(e_1)$ , and an edge $e_3 \in
G^1$ with $s(e_3) = r(e_2)$.  Continuing in this fashion we
will produce an infinite path $e_1e_2 \ldots $ with
$s(e_1)=v$.  Since $G$ is finite, the vertices $s(e_i)$ of
this path must eventually repeat.  Let $m$ be the smallest
natural number for which $s(e_m)=s(e_k)$ for some $1 \leq k
\leq m-1$.  Note that because $r(e_1) \neq s(e_1)$ we must
have $m \geq 3$.

Now $e_ke_{k+1} \ldots e_{n-1}$ will be a loop, and since $G$
satisfies Condition (L), there exists an exit for this loop. 
Thus for some $k \leq l \leq n-1$ there exists $f \in G^1$ such
that $r(f)=s(e_l)$ and $f \neq e_l$.  For each $w \in G^0$
define 
$$ a_w := \begin{cases} 2 & \text{if $w \in \{ s(e_i)
\}_{i=2}^l$} \\ 1 & \text{otherwise} \end{cases}$$
Note that $\{ s(e_i) \}_{i=2}^l$ may be empty.  This will
occur if and only if $l=1$.  Now let $n:=(a_w) \in \prod_{G^0}
\N$.  To see that $((A_G-I)n)(v) \geq 1$, note that $a_v =1$,
and consider four cases.

\noindent Case 1: $l=1$ and $r(f)=r(e_1)$.  Since $r(e_1) \neq
v$ we have that 
$$((A_G-I)n)(v) \geq a_v(A_G(v,v)-1) + a_{r(e_1)}A_G(v,r(e_1))
\geq 1(-1) + 1(2) =1.$$

\noindent Case 2: $l=1$ and $r(f)=v$.  Then
$$((A_G-I)n)(v) \geq a_v(A_G(v,v)-1) + a_{r(e_1)}A_G(v,r(e_1))
\geq 1(1-1) + 1(1) =1.$$

\noindent Case 3: $l=1$, $r(f) \neq r(e_1)$, and $r(f) \neq
v$.  Then
\begin{align}
((A_G-I)n)(v) \geq & \ a_v(A_G(v,v)-1)
+ a_{r(e_1)}A_G(v,r(e_1)) + a_{r(f)}A_G(v,r(f)) \notag \\ 
\geq & \ 1(-1) + 1(1) + 1(1) \notag \\
= & \ 1. \notag
\end{align}

\noindent Case 4: $l \geq 2$.  Then $a_{r(e_1)} =2$ and
$$((A_G-I)n)(v) \geq a_v(A_G(v,v)-1) + a_{r(e_1)}A_G(v,r(e_1))
\geq 1(-1) + 2(1) =1.$$

\noindent To see that $(A_G-I)n \in \prod_{G^0} \N$ let $w \in
G^0$ and consider the following three cases.

\noindent Case 1: $w=s(e_l)$ and $r(e_l)=r(f)$.  Then $a_w=2$
and we have 
$$((A_G-I)n)(w) \geq a_w(A_G(w,w)-1) + a_{r(e_l)}A_G(w,r(e_l))
\geq 2(-1) + 1(2) =0.$$

\noindent Case 2: $w=s(e_l)$ and $r(e_l) \neq r(f)$.  Then 
\begin{align}
((A_G-I)n)(w) \geq & \ a_w(A_G(w,w)-1) +
a_{r(e_l)}A_G(w,r(e_l)) + a_{r(f)}A_G(w,r(f)) \notag \\
\geq & \ 2(-1)+1(1)+1(1) \notag \\
= & \ 0. \notag
\end{align}

\noindent Case 3: $w \neq s(e_l)$.  Then either $w \in \{
s(e_i) \}_{i=2}^{l-1}$ or $a_w=1$.  In either case there
exists an edge $e$ with $s(e)=w$ and $a_{r(e)} \geq a_w$.  Thus
$$((A_G-I)n)(w) \geq a_w(A_G(w,w)-1) + a_{r(e)}A_G(w,r(e))
\geq -a_w+a_{r(e)} \geq 0$$
and $(A_G-I)n \in \prod_{G^0} \N$.  \end{proof}

\begin{theorem}   Let $G$ be a finite graph with no sinks
that satisfies Condition~(L).  For any $[x] \in \coker(A_G-I)$
there exists an essential 1-sink extension $E$ of $G$
such that $[\wvec_E] = [x]$ in $\coker(A_G-I)$.
\label{chooseonbesinkrep}
\end{theorem}
\begin{proof}  For each $v \in G^0$ we may use Lemma
\ref{finitevectorexists} to obtain a vector
$n_v \in \prod_{G^0} \N$ such that $(A_G-I)n_v \in
\prod_{G^0} \N$ and $((A_G-I)n_v)(v) \geq 1$.  Now write $x$ in
the form $x = \sum_{v \in G^0} a_v \delta_v$.  Let $n : =
\sum_{v \in G^0} ( |a_v| + 1 ) n_v$.  Then by linearity,
$x+(A_G-I)n \in \prod_{G^0} \N$ and $x+(A_G-I)n \neq 0$.  Let
$E$ be a 1-sink extension of $G$ with sink $v_0$ and $\W$ vector
equal to $x+(A_G-I)n$.  Then $[\wvec_E] = [x+(A_G-I)n] = [x]$ in
$\coker(A_G-I)$.  Furthermore, since $\wvec_E(v) \geq 1$ for
all $v \in G^0$ it follows that $G^0 \geq v_0$
and $E$ is an essential 1-sink extension. \end{proof}

This result shows that if $G$ is a finite graph with no sinks
that satisfies Condition~(L), then for any element in
$\ext (C^*(G))$ we may choose a representative that is the
extension associated to an essential 1-sink extension $E$ of
$G$.  Furthermore, since the $\W$ map is an isomorphism we see
that if $E_1$ and $E_2$ are essential 1-sink extensions that
are representatives for $\tau_1, \tau_2 \in \ext(C^*(G))$,
then the essential 1-sink extension with $\W$ vector equal to
$\wvec_{E_1} + \wvec_{E_2}$ will be a representative of $\tau_1
+ \tau_2$.  Hence we have a way of choosing representatives of
the classes in $\ext$ that have a nice visual interpretation
and for which we can easily compute their sum.

\section{Semiprojectivity of graph algebras}
\label{semiprojectivity}

In 1983 Effros and Kaminker \cite{EK} began the development of
a shape theory for $C^*$-algebras that generalized the
topological theory.  In their work they looked at
$C^*$-algebras with a property that they called
semiprojectivity.  These semiprojective
$C^*$-algebras are the noncommutative analogues of
absolute neighborhood retracts.  In 1985 Blackadar generalized
many of these results \cite{Bla2}, but because he wished
to apply shape theory to $C^*$-algebras not included in
\cite{EK} and because the theory in \cite{EK} was not a direct
noncommutative generalization, Blackadar gave a new definition
of semiprojectivity.  Blackadar's definition is more
restrictive than that in \cite{EK}.
\begin{definition}[Blackadar]  A separable $C^*$-algebra $A$
is \emph{semiprojective} if for any $C^*$-algebra $B$, any
increasing sequence $\{ J_n \}_{n=1}^\infty$ of (closed
two-sided) ideals, and any $*$-homomorphism $\phi : A
\rightarrow B / J$, where $J := \overline{\bigcup_{n=1}^\infty
J_n}$, there is an $n$ and a $*$-homomorphism $\psi : A
\rightarrow B / J_n$ such that $$ \xymatrix{A \ar[r]^-{\psi}
\ar[dr]_{\phi} &
B / J_n \ar[d]^{\pi} \\
 & B / J \\ } $$ where $\pi : B / J_n \rightarrow B / J$ is the
natural quotient map.
\end{definition}

In \cite{Bla2} it was shown that the Cuntz-Krieger algebras are
semiprojective, and more recently Blackadar has announced a
proof that $\mathcal{O}_\infty$ is semiprojective. 
Based on the proof for $\mathcal{O}_\infty$ Szyma\'nski has
proven in \cite{Szy} that if $E$ is a transitive graph with
finitely many vertices (but a possibly infinite number of
edges), then $C^*(E)$ is semiprojective.   

We now give an example of a row-finite transitive graph $G$
with an infinite number of vertices and with the property that
$C^*(G)$ is not semiprojective.  We use the fact that the
$\W$ map of \S \ref{Wmapsection} is an isomorphism in order to
prove that $C^*(G)$ is not semiprojective.

If $G$ is a graph, then by \emph{adding a sink at $v \in G^0$}
we shall mean adding a single vertex $v_0$ to $G^0$ and a single
edge $e$ to $G^1$ going from $v$ to $v_0$.  More formally, if
$G$ is a graph, then we form the graph $F$ defined by $F^0 :=
G^0 \cup \{ v_0 \}$, $F^1 := G^1 \cup \{ e \}$, and we
extend $r$ and $s$ to $F^1$ by defining and $r(e) = v_0$ and
$s(e)=v$.

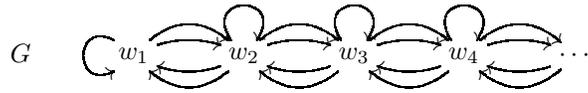
\begin{example}
\label{nonsemiprojex}
$$\begin{matrix}

G & & & 

\xymatrix{
w_1 \ar@(ul,dl) \ar@/^/[r] \ar@/^1pc/[r]  &
w_2 \ar@(ul,ur)  \ar@/^/[r] \ar@/^1pc/[r]  \ar@/^/[l]
\ar@/^1pc/[l] & w_3
\ar@(ul,ur)  \ar@/^/[r] \ar@/^1pc/[r]  \ar@/^/[l]
\ar@/^1pc/[l] & w_4 \ar@(ul,ur)  \ar@/^/[r] \ar@/^1pc/[r] 
\ar@/^/[l]
\ar@/^1pc/[l] & \cdots  \ar@/^/[l]
\ar@/^1pc/[l] \\
} \\ \end{matrix}$$

\vspace{.2in}

\noindent If $G$ is the above graph, then note that $G$
is transitive, row-finite, and has no sinks.
\end{example}

\begin{theorem}  If $G$ is the graph in
Example~\ref{nonsemiprojex}, then $C^*(G)$ is not
semiprojective.
\end{theorem}

\begin{proof}  For each $i \in \N$ let $E_i$ be the graph formed
by adding a sink to $G$ at $w_i$, and let $F_i$ be the graph
formed by adding a sink to each vertex in $\{ w_i, w_{i+1},
\ldots \}$.  In each case we shall let $v_i$ denote the sink
that is added at $w_i$.  As examples we draw $E_3$ and $F_3$:

$$\begin{matrix} E_3 & & &

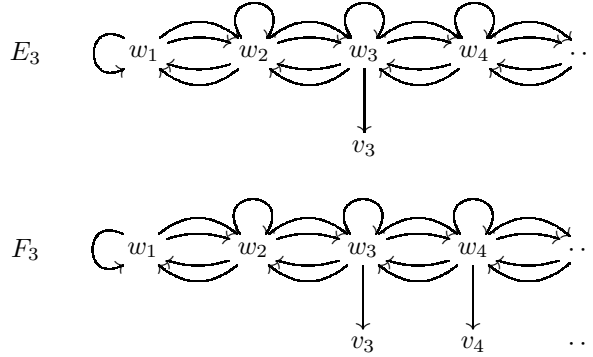
\xymatrix{
w_1 \ar@(ul,dl) \ar@/^/[r] \ar@/^1pc/[r]  &
w_2 \ar@(ul,ur)  \ar@/^/[r] \ar@/^1pc/[r]  \ar@/^/[l]
\ar@/^1pc/[l] & w_3 \ar[d]
\ar@(ul,ur)  \ar@/^/[r] \ar@/^1pc/[r]  \ar@/^/[l]
\ar@/^1pc/[l] & w_4 \ar@(ul,ur)  \ar@/^/[r] \ar@/^1pc/[r] 
\ar@/^/[l]
\ar@/^1pc/[l] & \cdots  \ar@/^/[l]
\ar@/^1pc/[l] \\
 & & v_3 & & \\
} \\ \end{matrix}$$

\vspace{.3in}

$$\begin{matrix} F_3 & & &

\xymatrix{
w_1 \ar@(ul,dl) \ar@/^/[r] \ar@/^1pc/[r]  &
w_2 \ar@(ul,ur)  \ar@/^/[r] \ar@/^1pc/[r]  \ar@/^/[l]
\ar@/^1pc/[l] & w_3 \ar[d]
\ar@(ul,ur)  \ar@/^/[r] \ar@/^1pc/[r]  \ar@/^/[l]
\ar@/^1pc/[l] & w_4 \ar[d] \ar@(ul,ur)  \ar@/^/[r]
\ar@/^1pc/[r] 
\ar@/^/[l]
\ar@/^1pc/[l] & \cdots  \ar@/^/[l]
\ar@/^1pc/[l] \\
 & & v_3 & v_4 & \cdots \\
}  \\ \end{matrix}$$

We shall now assume that $C^*(G)$ is semiprojective and arrive
at a contradiction.  Let $B : = C^*(F_1)$ and for each $n
\in \N$ let $H_n := \{ v_1, v_2, \ldots , v_n \}$.  Also let
$H_\infty := \{ v_1, v_2, \ldots  \}$.  Set $J_n := I_{H_n}$. 
Then $\{ J_n \}_{n=1}^\infty$ is an increasing sequence of
ideals and $J := \bigcup_{n=1}^\infty J_n = I_{H_\infty}$. 
Now $B / J = C^*(F_1) / I_{H_\infty} \cong C^*(G)$ and for
each $n
\in \N$, $B / J_n
\cong C^*(F_{n+1})$ by
\cite[Theorem 4.1]{BPRS}.  Thus if we identify $C^*(G)$ and $B
/ J$, then by semiprojectivity there exists a homomorphism
$\psi : C^*(G) \rightarrow B / J_n$ for some $n$ 
$$ \xymatrix{C^*(G) \ar[r]^-{\psi}
\ar[dr]_{\text{id}} &
B / J_n \cong C^*(F_{n+1}) \ar[d]^{\pi} \\
 & B / J \cong C^*(G) \\ } $$ such that $\pi \circ \psi =
\text{id}$.  Note that the projection $\pi : B /J_n
\rightarrow B/J$ is just the projection $\pi : C^*(F_{n+1})
\rightarrow C^*(F_{n+1}) / I_{ \{ v_{n+1},v_{n+2},\ldots \} }
\cong C^*(G)$.

Now if we let $\{s_e, p_v
\}$ be the canonical Cuntz-Krieger $F_{n+1}$-family in
$C^*(F_{n+1})$ and let $\{t_e, q_v \}$ be the canonical
Cuntz-Krieger $E_{n+1}$-family in $C^*(E_{n+1})$, then by the
universal property of $C^*(F_{n+1})$ there exists a
homomorphism $\rho : C^*(F_{n+1}) \rightarrow C^*(E_{n+1})$
such that $$\rho(s_e) = \begin{cases} t_e & \text{ if $e \in
E_{n+1}^1$} \\ 0 & \text{ if $e \in F_{n+1}^1 \backslash
E_{n+1}^1$ } \end{cases} \hspace{.2in} \text{and} \hspace{.2in}
\rho(p_v) = \begin{cases} q_v & \text{ if $v \in E_{n+1}^0$} \\
0 & \text{ if $v \in F_{n+1}^0 \backslash E_{n+1}^0$. }
\end{cases}$$ Since $E_{n+1}$ is a 1-sink extension of $G$, we
have the usual projection $\pi_{E_{n+1}} : C^*(E_{n+1})
\rightarrow C^*(G)$.  One can then check that the diagram
$$ \xymatrix{C^*(F_{n+1}) \ar[rr]^{\rho}
\ar[dr]_{\pi} & &
C^*(E_{n+1}) \ar[dl]^{\pi_{E_{n+1}}} \\
 & C^*(G) & \\ } $$ commutes simply by checking that
$\pi_{E_{n+1}} \circ \rho$ and
$\pi$ agree on generators.  This, combined with the fact that
$\pi \circ \psi = \text{id}$ on $C^*(G)$, implies that
$\pi_{E_{n+1}} \circ \rho \circ \psi = \text{id}$.  Hence the
short exact sequence $$\xymatrix{ 0 \ar[r] & I_{v_{n+1}}
\ar[r] & C^*(E_{n+1}) \ar[r]^{\pi_{E_{n+1}}} & C^*(G) \ar[r]
\ar@(ul,ur)[l]_{\rho \circ \psi} & 0}$$ is split exact. 
Therefore this extension is degenerate.  Since $I_{v_{n+1}}
\cong \K$ by \cite[Corollary 2.2]{KPR} we have that this
extension is in the zero class in $\ext(C^*(G))$.

However, the $\W$ vector of $E_{n+1}$ is
$\wvec_{E_{n+1}} =
\delta_{w_{n+1}}$.  Since $$ A_G-I = \left( \begin{smallmatrix}
0 & 2 & 0 & 0 & \\
2 & 0 & 2 & 0 & \cdots \\
0 & 2 & 0 & 2 & \\
0 & 0 & 2 & 0 & \\
  & \vdots  &  &  & \ddots
\end{smallmatrix} \right)$$ we see that every vector in the
image of $A_G-I$ has entries that are multiples of $2$.  Thus
$\delta_{w_{n+1}} \notin \im (A_G-I)$, and $[\wvec_{E_{n+1}}]$
is not zero in $\coker (A_G-I)$.  But then Proposition
\ref{wclassequalswvec} and Theorem \ref{wisaniso} imply that
the extension associated to $C^*(E_{n+1})$ is not equal to
zero in $\ext(C^*(G))$.  This provides the contradiction, and
hence $C^*(G)$ cannot be semiprojective. \end{proof}

\begin{remark}
After the completion of this work, Spielberg
proved in
\cite{Spi} that all classifiable, simple, separable, purely
infinite $C^*$-algebras having finitely generated $K$-theory
and free $K_1$-group are semiprojective
\cite[Theorem~3.12]{Spi}.  This was accomplished by realizing
these
$C^*$-algebras as graph algebras of transitive
graphs.  It also implies that if $G$ is a transitive graph that
is not a single loop, and if $C^*(G)$ has finitely generated
$K$-theory and free $K_1$-group, then $C^*(G)$ is
semiprojective.  We mention that the $C^*$-algebra associated
to the graph in Example~\ref{nonsemiprojex} does not have
finitely generated
$K$-theory.
\end{remark}


\begin{thebibliography}{99}


\bibitem{Arv} 
W. Arveson, \emph{Notes on extensions of $C^*$-algebras}, Duke
Math. J. \textbf{44} (1977), 329--355.

\bibitem{BPRS}
T. Bates, D. Pask, I. Raeburn, and W. Szyma\'nski,
\emph{The $C^*$-algebras of row-finite graphs}, New York J.
Math. \textbf{6} (2000), 307--324.

\bibitem{Bla}
B. Blackadar, K-theory for Operator Algebras, second ed.,
Cambridge University Press., Cambridge 1998. 

\bibitem{Bla2}
B. Blackadar, \emph{Shape theory for $C^*$-algebras}, Math.
Scand. \textbf{56} (1985), 249--275  

\bibitem{CK} 
J. Cuntz and W. Krieger, \emph{A class of $C^*$-algebras and
topological markov chains}, Invent. Math. \textbf{56} (1980),
no. 3, 251--268

\bibitem{Dav} 
K. Davidson, $C^*$-algebras by Example, Fields Institute
Monographs, vol. 6, Amer. Math. Soc., Providence, 1996.

\bibitem{EK}
E. G. Effros and J. Kaminker, \emph{Homotopy continuity and
shape theory for $C^*$-algebras}, in Geometric Methods in
Operator Algebras (Kyoto 1983), pp. 152--180, Pitman Res. Notes
Math. \textbf{123}, Longman Sci. Tech., Harlow, 1986. 

\bibitem{JT}
K. K. Jensen and K. Thomsen, Elements of
$KK$-theory, Birkh\"auser, Boston, Massachusetts, 1991.

\bibitem{KPR}
A. Kumjian, D. Pask, and I. Raeburn, \emph{Cuntz-Krieger
algebras of directed graphs}, Pacific J. Math. \textbf{184}
(1998), 161--174.

\bibitem{KPRR}
A. Kumjian, D. Pask, I. Raeburn and J. Renault,
\emph{Graphs, groupoids, and Cuntz-Krieger algebras}, J. Funct.
Anal. \textbf{144} (1997), 505--541.

\bibitem{MRS}
M. H. Mann, I. Raeburn, C. E. Sutherland,
\emph{Representations of finite groups and Cuntz-Krieger
algebras}, Bull. Austral. Math. Soc.  \textbf{46} (1992),
225--243. 

\bibitem{RS}
I. Raeburn and W. Szyma\'nski,
\emph{Cuntz-Krieger algebras of infinite graphs and matrices},
Proc. Amer. Math. Soc. \textbf{129} (2000), 2319--2327.

\bibitem{RTW}
I. Raeburn, M. Tomforde, D. Williams,
\emph{Classification theorems for the $C^*$-algebras of graphs
with sinks}, preprint (2000).

\bibitem{RW}
I. Raeburn and D. Williams, Morita Equivalence and
Continuous-Trace $C^*$-algebras, Math. Surveys \& Monographs,
vol. 60, Amer. Math. Soc., Providence, 1998.

\bibitem{Spi}
J.~Spielberg, \emph{Semiprojectivity for certain purely infinite
$C^*$-algebras}, preprint (2001).

\bibitem{Szy}
W. Szyma\'nski, \emph{On semiprojectivity of $C^*$-algebras of
directed graphs}, Proc. Amer. Math Soc., \textbf{130} (2002),
1391--1399.

\bibitem{thesis}
M.~Tomforde, \emph{Extensions of graph $C^*$-algebras}, Ph.D.
thesis, Dartmouth College, 2002.

\bibitem{Voi} 
D. Voiculescu, \emph{A non-commutative Weyl-von Neumannn
theorem}, Rev. Roumaine Math. Pures Appl. \textbf{21} (1976),
97--113.  

\bibitem{WO}
N. E. Wegge-Olsen, $K$-theory and $C^*$-algebras, Oxford
University Press, Oxford, 1993.  

\bibitem{Wat}
Y. Watatani, \emph{Graph theory for $C^*$-algebras}, in
Operator Algebras and Their Applications (R.V. Kadison, ed.),
Prpc. Symp. Pure Math., vol. 38, part~1, Amer. Math. Soc.,
Providence, 1982, pages 195--197. 

\end{thebibliography}
\end{document}